\documentclass[10pt]{amsart}
\usepackage{graphicx,amssymb,amsfonts,amsmath,amsthm,newlfont}
\usepackage{epsfig}
\newtheorem{thm}{Theorem}
\newtheorem{cor}[thm]{Corollary}
\newtheorem{lem}[thm]{Lemma}
\newtheorem{prop}[thm]{Proposition}

\theoremstyle{definition}
\newtheorem{defn}[thm]{Definition}
\newtheorem{notation}[thm]{Notation}

\newtheorem{example}[thm]{Example}

\theoremstyle{remark}
\newtheorem{rem}[thm]{Remark}

\font \rus= wncyr10

\newcommand{\Z}{\mathbb Z}
\newcommand{\C}{\mathbb C}
\newcommand{\Pro}{\mathbb P}

\newcommand{\R}{\mathbb R}
\newcommand{\N}{\mathbb N}
\newcommand{\Q}{\mathbb Q}

\newcommand{\Reg}{\mathrm{Reg}}

\newcommand{\q}{/\!\!/}

\newcommand{\x}{\mathsf{x}}
\newcommand{\ao}{\mathsf{a}}
\newcommand{\MZV}{\mathcal{Z}}
\newcommand{\RReg}{\mathrm{RReg}}
\newcommand{\V}{\mathcal V}

\newcommand{\To}{\longrightarrow}
\newcommand{\sha}{\, \hbox{\rus x} \,}

\newcommand{\Mod}{\mathfrak{M}}

\newcommand{\Sym}{\mathfrak{S}}

\newcommand{\Or}{\mathcal{O}}

\newcommand{\Real}{\mathrm{Re\,}}

\newcommand{\y}{ \mathsf{y}}

\newcommand{\Li}{\mathrm{Li}}

\begin{document}
\title{The Massless higher-loop two-point function}
\author{Francis Brown, CNRS, IMJ}
\date{10 April 2008} \maketitle
\begin{abstract} We introduce a new method for computing massless Feynman
integrals analytically in parametric form. 
An analysis of the method yields a criterion for a primitive Feynman
graph $G$ to evaluate to multiple zeta values. The criterion depends
only on the topology of $G$, and can be checked algorithmically. As
a corollary, we reprove the result, due to Bierenbaum and Weinzierl,
that the massless 2-loop 2-point function is expressible in terms of
multiple zeta values, and generalize this to  the 3, 4, and 5-loop cases. We find that the coefficients in the Taylor expansion
of planar  graphs in this  range  evaluate to multiple zeta values, but the non-planar  graphs  with crossing number 1 may evaluate
 to multiple sums  with $6^\mathrm{th}$ roots of unity.
Our method  fails for the five loop graphs with crossing number 2 obtained by breaking open the  bipartite   graph  $K_{3,4}$ at one edge.
%
\end{abstract}

\section{Introduction}
Let $n_1,\ldots, n_r\in \N$ and suppose that $n_r\geq 2$. The
multiple zeta value is the real number defined by the convergent
nested sum:
$$\zeta(n_1,\ldots, n_r) =\sum_{0< k_1< k_2<\ldots<k_r } {1 \over k_1^{n_1}\ldots k_r^{n_r}} \ .$$
An important question in perturbative quantum field theory is
whether multiple zeta values, or some larger set of periods, suffice
to evaluate a given class of Feynman integrals. Moreover, it is
crucial for applications to find efficient methods for evaluating
such Feynman integrals analytically.

In this paper, we  shall  
 consider 
massless Feynman integrals in $\phi^4_4$ theory with
the propagators raised to arbitrary powers, and a single non-zero momentum. The simplest case is the
two-loop example, pictured in figure 1 on the left. Let $a_i$ be a
positive real number corresponding to each edge $1\leq i\leq 5$.
\begin{figure}[h!]\label{figureThreeLoop}
  \begin{center}
    \epsfxsize=10.0cm \epsfbox{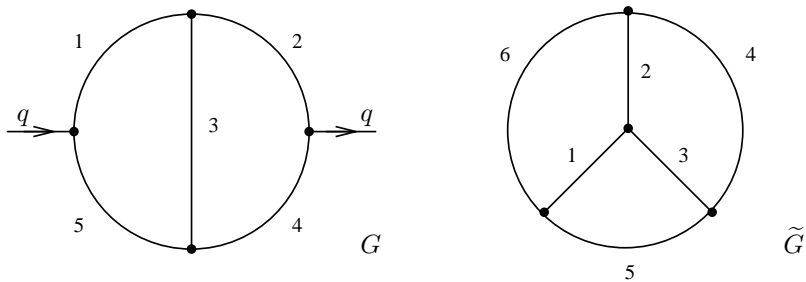}
\put(-280,60){$q$} \put(-150,60){$q$}
\put(-150,10){$G$} \put(10,10){$\widetilde{G}$}
 \caption{Left: The two-point two-loop massless diagram. Closing up its external legs gives the
wheel with three spokes  (right).}
  \end{center}
\end{figure}
The corresponding Feynman integral is:
\begin{equation} \label{2P2LFeynman}
\int \!\!\int d^D k_1 d^D k_2 \Big({1 \over k_1^2}\Big)^{a_1}
\Big({1\over k_2^2}\Big)^{a_2}
 \Big({1\over (k_1-k_2)^2}\Big)^{a_3} \Big({1 \over (k_2-q)^2}\Big)^{a_4}\Big({1 \over (k_1-q)^2}\Big)^{\alpha_5}\ ,
\end{equation}
where $D$ is the number of dimensions, and $q$ is the momentum
entering on the left. Now suppose that
 $a_i=1 + n_i \,\varepsilon$, where $n_i$ are positive integers,
for $1\leq i\leq 5$, and $\varepsilon$ is the parameter in
dimensional regularization, i.e.,  $D=4-2\,\varepsilon$. The
 problem of calculating the coefficients in the Taylor
expansion of $(\ref{2P2LFeynman})$ with respect to $\varepsilon$, is
important in three and four loop calculations, and has a history
spanning approximately twenty-five years (to which we refer to
\cite{B-W} for an account). In particular, it had been conjectured
for a long time that every coefficient is a rational linear
combination of multiple zeta values. This question was finally
settled in the affirmative in \cite{B-W}, using Mellin-Barnes
techniques.

Until now, however, there seemed to be a lack of systematic methods
for computing a range of Feynman integrals analytically at higher
loop orders. In this paper, we introduce a new method, using
iterated integration with polylogarithms, which was initiated in
\cite{Br3} to compute the periods of moduli spaces of curves of
genus 0. Using this, we reprove the fact that the Taylor expansion
of the massless two-loop two-point integral $(\ref{2P2LFeynman})$
evaluates to multiple zeta values, and extend the result to higher
loop orders. Our method also yields results for certain examples of
massive Feynman diagrams, but in the present paper we only consider
 massless cases.



\subsection{Results}
We consider three and higher-loop integrals with exactly one non-zero momentum and arbitrary powers of the propagators, generalizing the
integral $(\ref{2P2LFeynman})$. Let $G$ be a graph in $\phi^4_4$
with two external legs, and a single momentum $q$. By Hopf algebra
arguments, we can assume that $G$ is primitive, in the following
sense. Let $\widetilde{G}$ denote the graph formed by closing the
two external legs of $G$, which now has no external edges, but gains
an extra loop (see figure 1).  We will suppose that
$\widetilde{G}$ is primitive divergent in the sense of \cite{B-E-K},
i.e., $\widetilde{G}$ contains no strict divergent subgraphs and
satisfies:
$$\hbox{\#edges of } \widetilde{G} = 2\times  \hbox{\#loops of } \widetilde{G}\ .$$
In this case, we say that $G$ is broken primitive divergent
(bpd), and every bpd graph with a fixed number of loops can be
obtained by taking the set of all primitive divergent graphs with
one more loop, and breaking them open along every edge.

Now  suppose that $G$ is bpd with $h$ loops and $L$ internal edges. It
follows from the above that $L=2h+1$.  For each internal edge $i$,
let $a_i= 1 + n_i\,\varepsilon$, where $n_i$ is a positive integer,
and $D=4-2\,\varepsilon$, and consider the massless Feynman integral
\begin{equation} \label{intropreFeyn}
\int d^Dk_1\ldots \int  d^Dk_L  \prod_{i=1}^L \,\Big( {1\over
r_i}\Big)^{a_i}\ ,
 \end{equation}
where $k_i$ is a momentum flowing through the
$i^{\mathrm{th}}$ edge, and $r_i$ is the propagator corresponding to
the $i^{\mathrm{th}}$ edge. The domain of integration is given by the conservation of momentum at each vertex,
and there is a single external momentum entering $G$, denoted $q$.
It is easy to show, using the Schwinger
trick, that $(\ref{intropreFeyn})$ is proportional to a certain
power of $q^2$ \cite{IZ}. 
The
constant of proportionality  is given by  a product of explicit
gamma factors and powers of $\pi$, with the parametric
 integral:
\begin{equation} \label{introgeneralFeyn}
I(G) = \int_{\alpha_{\lambda}=1} {\prod_{i=1}^{L+1} \alpha_i^{a_i-1}
\over U_{\widetilde{G}}^{D/2} }\, d\alpha_1\ldots
\widehat{d\alpha}_\lambda \ldots d\alpha_{L+1}\ ,
 \end{equation}
where
 $U_{\widetilde{G}}(\alpha_1,\ldots,
\alpha_{L+1})$ is the graph polynomial (or Kirchoff polynomial) of
$\widetilde{G}$, $\alpha_1,\ldots, \alpha_{L+1}$ are  Schwinger
parameters for the $L+1$ edges in $\widetilde{G}$, $\lambda$ is any
index between $1$ and $L+1$, and
  $a_{L+1}$ is a certain linear combination of
$a_1,\ldots, a_L$. The problem of computing the Taylor expansion of
$(\ref{intropreFeyn})$ reduces to computing the expansion of
$(\ref{introgeneralFeyn})$.  For example, if $\widetilde{G}$ is the
wheel with 3 spokes, we have:
\begin{eqnarray} \label{Ugtilde3spokes}
&U_{\widetilde{G}}=
&\!\!\!\alpha_1\alpha_2\alpha_6\!+\!\alpha_1\alpha_4\alpha_6\!+\!\alpha_2\alpha_5\alpha_6\!+\!\alpha_4\alpha_5\alpha_6
\!+\!
\alpha_1\alpha_3\alpha_6\!+\!\alpha_2\alpha_3\alpha_6\!+\!\alpha_3\alpha_4\alpha_6\!+\!\alpha_3\alpha_5\alpha_6  \nonumber \\
& \quad +
&\!\!\!\alpha_1\alpha_3\alpha_4\!+\!\alpha_1\alpha_3\alpha_5\!+\!\alpha_2\alpha_3\alpha_4\!+\!\alpha_2\alpha_3\alpha_5\!+\!
\alpha_2\alpha_4\alpha_5\!+\!\alpha_1\alpha_4\alpha_5\!+\!\alpha_1\alpha_2\alpha_5\!+\!\alpha_1\alpha_2\alpha_4
 \nonumber \ ,
\end{eqnarray}
and the coefficients in the Taylor expansion of
$(\ref{introgeneralFeyn})$ are given by the period integrals:
\begin{eqnarray}\label{Intrologint}
\int_{\alpha_{\lambda}=1}  {\log(\alpha_1)^{m_1}\ldots
\log(\alpha_{L+1})^{m_{L+1}}   \log(U_{\widetilde{G}})^n  \over
U_{\widetilde{G}}^{2} } \, d\alpha_1\ldots \widehat{d\alpha}_\lambda
\ldots d\alpha_{L+1}\ ,
\end{eqnarray}
where $m_1,\ldots, m_{L+1}, n$ are arbitrary positive integers.

 We know state some results on
the transcendental nature of the
 coefficients  in the Taylor expansion of
$I(G)$ with respect to $\varepsilon$, for all bpd graphs $G$ up to
five loops.
For every graph $G$ for which a theorem is  stated below, there is
also a corresponding   algorithm for computing the coefficients of
the Taylor expansion of $I(G)$ by integrating inside  a
predetermined algebra of polylogarithms.


\subsubsection{Three loops} There is
exactly one primitive divergent graph with four loops  \cite{KY},
namely the wheel with four spokes, or cross-hairs diagram, pictured
below (left). Breaking it apart along each edge gives rise to
exactly two topologically distinct  bpd graphs with three loops
(right).

\begin{figure}[h!]
  \begin{center}
    \epsfxsize=13.0cm \epsfbox{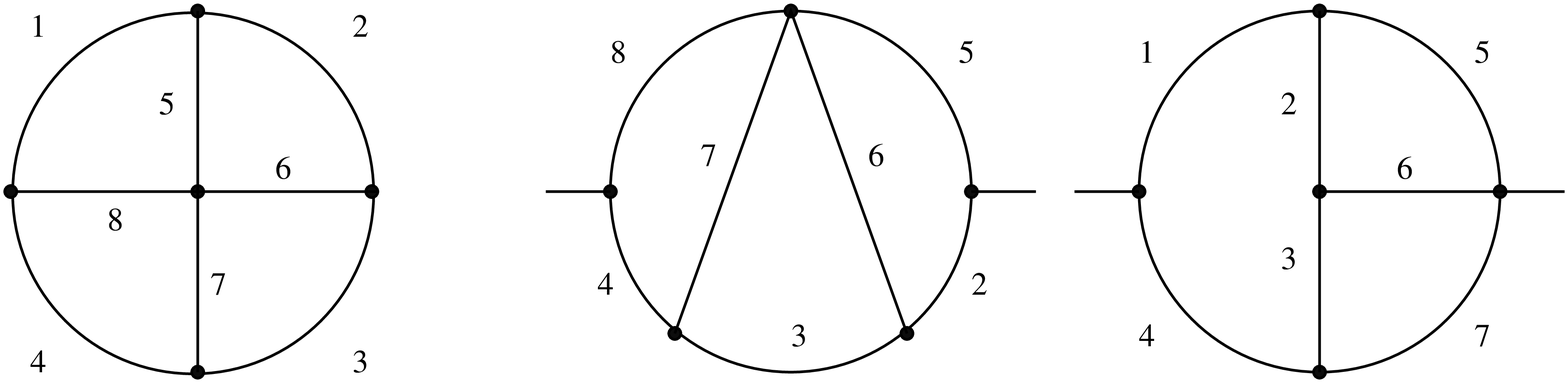}
  \label{figure1}
 \caption{Left: The cross-hairs diagram is the unique primitive divergent graph with four loops. Breaking it at the edges $1$, $8$ gives the graphs in the middle, and on the right, respectively.
}
  \end{center}
\end{figure}

\begin{thm} Let $G$ be one  of the two bpd 3-loop graphs depicted in figure 2 (middle and right). Then every coefficient in the Taylor expansion
of $I(G)$ is a rational linear combination of multiple zeta values.
\end{thm}

\subsubsection{Four loop contributions} There are precisely three
5-loop diagrams, pictured in figure 3. The one on the left is
planar and two-vertex reducible, and will be denoted $5\,R$. The one
in the middle is planar, two-vertex irreducible, and will be denoted
$5\,P$,  and  the one on the  right is non-planar, and will be
denoted $5\,N$.
\begin{figure}[h!]
  \begin{center}
    \epsfxsize=12.0cm \epsfbox{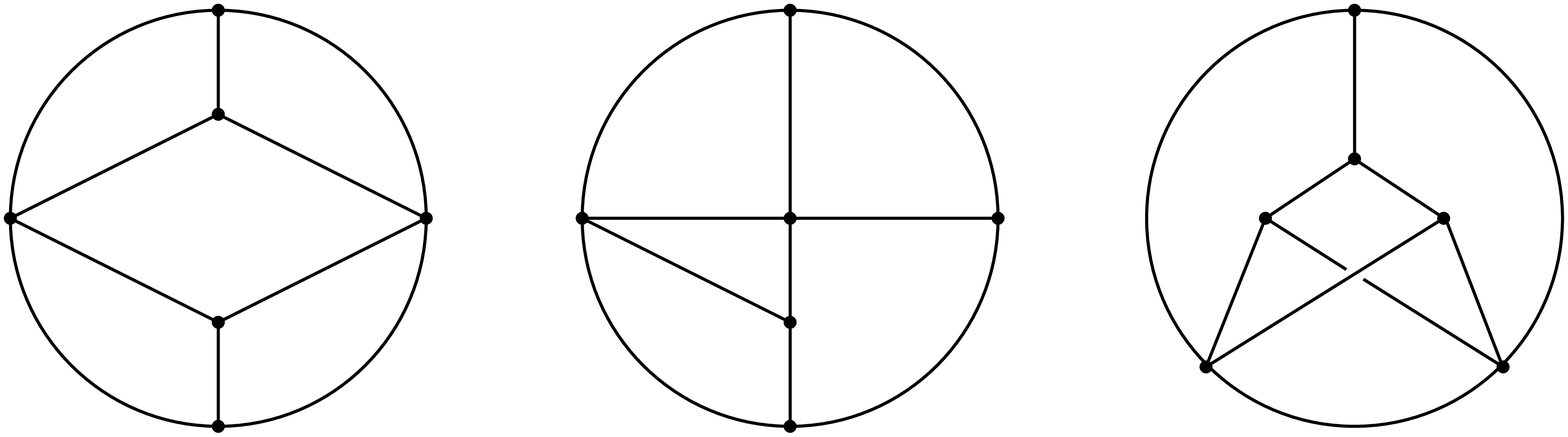}
  \label{figure2}
\put(-260,5){\small $5\,R$} \put(-135,5){\small $5\,P$}
\put(-10,5){\small $5\,N$}
 %
\caption{The three primitive-divergent 5-loop diagrams.}
  \end{center}
\end{figure}
It turns out that there are exactly six topologically distinct  ways
to break the graph $5\,P$ open at one edge, and exactly two ways to
break open $5\,R$, giving the  eight planar topologies depicted in
figure 4.
\begin{figure}[h!]
  \begin{center}
    \epsfxsize=12.0cm \epsfbox{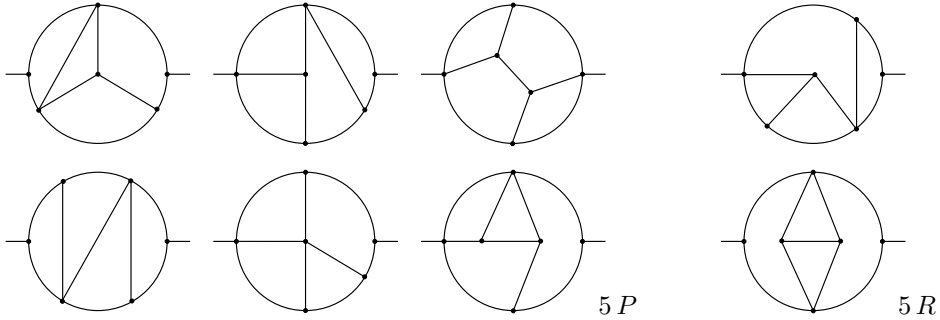}
  \label{figure2}
\put(-120,0){$5\,P$} \put(-10,0){ $5\,R$}
%
 \caption{The eight planar bpd topologies with four loops. 
}
  \end{center}
\end{figure}


\begin{thm}
Let $G$ be one of the eight bpd planar 4-loop graphs depicted in figure
4. Then every coefficient in the Taylor expansion of $I(G)$ is a
rational linear combination of multiple zeta values.\end{thm}

For the non-planar graphs, a new phenomenon occurs, and we must
introduce multiple zeta values at roots of unity. Let $n_1,\ldots,
n_r\in \N$, and consider the multiple polylogarithm function, first
introduced by Goncharov \cite{Go1}:
\begin{equation}\label{introMultPolydef}
\Li_{n_1,\ldots, n_r}(x_1,\ldots,x_r)= \sum_{0<k_1<\ldots<k_r}
{x_1^{k_1}\ldots x_r^{k_r}\over k_1^{n_1}\ldots k_r^{n_r}}\ .
\end{equation}
It converges absolutely for $|x_i|<1$ and extends to a multivalued
holomorphic function on an open subset of $\C^r$.
For $m\geq 1$, we define $\MZV^m$ to be the $\Q$-algebra generated
by the values of multiple polylogarithms at $m^{\mathrm{th}}$ roots
of unity:
$$\Li_{n_1,\ldots, n_r}(x_1,\ldots,x_r)\,\,\hbox{ such that } x_i^m=1 \hbox{ for } 1\leq i\leq r\ , \hbox{ and } (x_r,n_r)\neq (1,1)\ .$$
The condition that $x_r$ and $n_r$ are not simultaneously $1$ is to ensure  convergence.
The algebra $\MZV^1$ is the algebra of multiple
zeta values, and $\MZV^2$ is known as the algebra of alternating
multiple zeta values. We will call $\MZV^m$ the algebra of multiple
zeta values ramified at $m^{\mathrm{th}}$ roots of unity. Note that
$\MZV^a\subset \MZV^b$ if a divides b.

\begin{figure}[h!]
  \begin{center}
    \epsfxsize=10.0cm \epsfbox{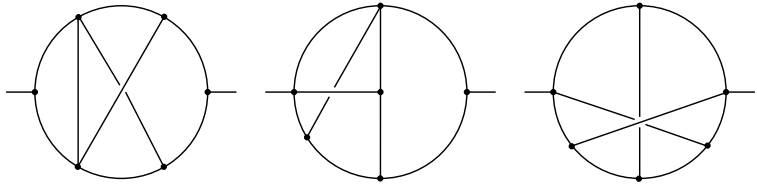}
  \label{figure2}
 \caption{The three non-planar bpd graphs on four loops obtained by breaking the graph $5\,N$ along an edge. The graph on the right has crossing number 1, but
has been  drawn with 2 crossings.}
  \end{center}
\end{figure}

\begin{thm}Let $G$ be one of the three non-planar bpd graphs with
4-loops as depicted in figure 5. Then every coefficient in the
Taylor expansion of $I(G)$ is a rational linear combination of
multiple
zeta values ramified at $6^{\mathrm{th}}$ roots of unity. 
\end{thm}

Note that the theorem gives an upper bound on the set of periods
which can occur ($\MZV^6$), and  it is an open question whether some
smaller algebra $\MZV^k$ for $k=1,2,\hbox{ or } 3$ suffices to
compute the Taylor expansions for the non-planar graphs, or whether
there occurs a term which genuinely involves  sixth roots of unity.

\subsubsection{Five-loop contributions}
 Karen Yeats has computed the 
 primitive divergent topologies up to seven loops \cite{KY}. There are nine at six loops, of which four  are planar, pictured below,
and  five which are non-planar (figures 7,8). Although the five non-planar graphs all have genus 1 (they can be drawn on a torus without any self-crossings),
it turns out that the invariant which succesfully predicts the outcome of our method for computing the periods at six loops is  not the genus but the crossing number.
The crossing number  of a graph $G$ is defined to be the minimial number of self-crossings over all planar representations of $G$ (which are allowed to have
curved edges). It is easy to determine the crossing number for the graphs we consider here, but seems to be a  difficult problem in general.
\begin{thm} Let $G$ be a graph obtained by breaking open $\widetilde{G}$ at one edge, where $\widetilde{G}$ is one of the eight primitive divergent graphs at six loops pictured in figures
6 and 7. Then the coefficients in the Taylor expansion of $I(G)$ are:
\begin{enumerate}
  \item multiple zeta values, if $\,\widetilde{G}$ is planar (has crossing number 0).
  \item multiple zeta values at $6^\mathrm{th}$ roots of unity, if $\widetilde{G}$ has crossing number 1\ .
\end{enumerate}
\end{thm}

The one remaining primitive  graph is the complete bipartite graph $K_{3,4}$, and has crossing number exactly 2, pictured in figure 8.
 Our method failed for this graph.

\begin{figure}[h!]
  \begin{center}
    \epsfxsize=12.0cm \epsfbox{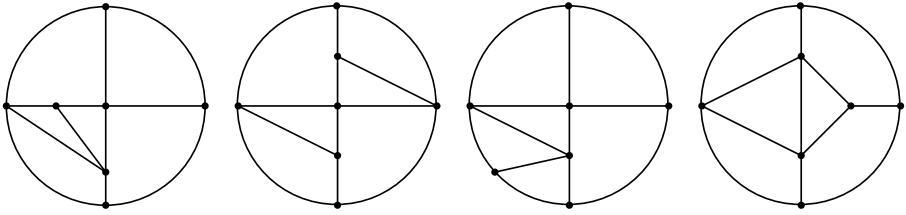}
  \label{figure2}
 \caption{The four planar primitive divergent graphs at six loops. }
  \end{center}
\end{figure}

\begin{figure}[h!]
  \begin{center}
    \epsfxsize=12.0cm \epsfbox{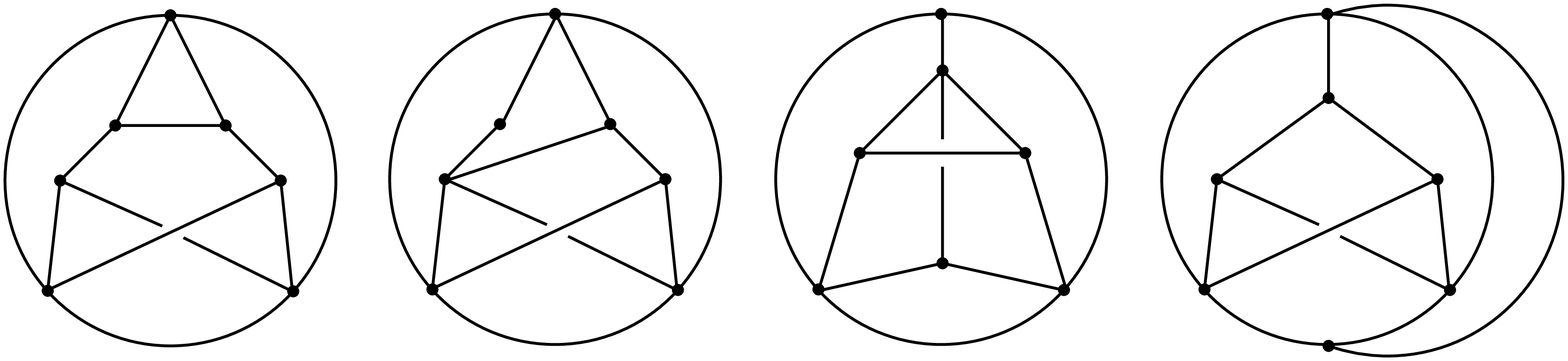}
  \label{figure2}
 \caption{The four non-planar primitive divergent graphs at six loops with crossing number 1.}
  \end{center}
\end{figure}

\begin{figure}[h!]
  \begin{center}
    \epsfxsize=3.0cm \epsfbox{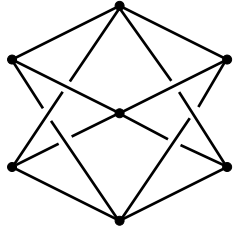}
  \label{figure2}
 \caption{The unique non-planar primitive divergent graph  at six loops with crossing number 2 (although six crossings are shown). It is the complete
bipartite graph $K_{3,4}$, and  is the first graph for which
our criterion fails.}
  \end{center}
\end{figure}

\subsection{Discussion}
In summary,  we found that for all Feynman graphs up to 6 loops, those  with crossing number 0 evaluate to multiple zeta values, and those with crossing number
 1 give  multiple polylogarithms evaluated  at $6^\mathrm{th}$ roots of unity. There is a single example of a graph ($K_{3,4}$) with crossing number 2, and our
method fails in this case.
The leading term in the Taylor expansion of $I(G)$ for $K_{3,4}$ is actually known to evaluate to a multiple zeta value, namely $\zeta(5,3)$, which, interestingly,  is the first
occurrence of an irreducible double sum in $\phi^4$ theory.

One possible reason for this is that our method concerns all terms in the Taylor expansion of the integral $I(G)$, with arbitrary powers of logarithms in the numerator. In the particular
case when  there are no logarithmic
terms in the numerator, the set of singularities of the integrand are slightly reduced, and it is not impossible that our method might work for $K_{3,4}$ with this restriction, although we have not checked this.
Another possible reason for this could be because our algorithm has some room for improvement, in at least two ways. First, the algorithm  involves the repeated factorization of polynomials derived from the graph polynomial.
In our computations, we only considered factorizations that occured over the field of rationals $\Q$, although it is conceivable that some of the  polynomials
 which occur are irreducible over $\Q$, but factorize over an algebraic extension of $\Q$.

 A second, and  more promising, possibility is to extend our algorithm to deal with
quadratic terms, and this will be explained in further detail in $\S\ref{sectExtensions}.$ It is possible, by extending our method in this way,
that
a larger class of graphs will become tractable. However, it is more likely than not that eventually one will find periods of motives which are not mixed Tate
({\it c.f.} \cite{BB}), and this will pose a genuine obstruction to the present method.
We expect that our method of polylogarithmic integration should also help to exhibit the first example of such a period in massless $\phi^4_4$, if and when
it occurs. The idea would be to
strip away  from a candidate Feynman integral terms which evaluate to multiple zeta values, until one is left with a totally irreducible period integral which is verifiably
not of mixed Tate type.  Thus, the eventual failure of our method  should at the same time exhibit the first non-MZV-type period, when it occurs.



\subsection{Plan of the paper}
The paper is divided into two halves. The first half (\S 2-4) gives
an overview of  our method, and can be read linearly, as
 a long introduction. In $\S2$ we briefly recall how to
rewrite the Feynman integrals $I(G)$  in Schwinger-parametric  form
using graph polynomials. 
In $\S3$, we outline the main idea of our method, and in $\S4$, we
translate this into an elementary reduction algorithm on graph
polynomials. The main theorem $\ref{thmMAINCRITERION}$ in
$\S\ref{sectMainthm}$ states a sufficient condition for a bpd graph
$G$ to evaluate to multiple zeta values in its Taylor expansion.

The purpose of the second half is to give a complete worked example
of our method, in the case of the wheel with three spokes. 
 In particular, in \S7, we give a new proof of the well-known
result that the leading term of $I(G)$ is $6\,\zeta(3)$. This is, to
our knowledge, the first time that such a computation has appeared
in print  using a parametric representation for the Feynman
integral. The purpose of section 5 is to provide sufficiently many
details to understand the intricacies of the method in general, and
 section 6 provides  worked examples of taking primitives and limits of
polylogarithms in two variables, which are then applied directly in
the $6\,\zeta(3)$ computation. It is perhaps advisable, after having
read $\S2-4$, to refer directly to the example in $\S7$, and then
read $\S5$ and $\S6$ bearing the example in mind.
\\

The proofs of our results  will be written up in full detail in
\cite{Br4}, including a study  of the underlying algebraic geometry,
which is completely absent from this paper, and a generalization of the above results to some infinite families of graphs. Much of the background
on iterated integrals, and also our method of integration with
polylogarithms, is explained in detail in \cite{Br2,Br3} in the case
of moduli spaces of genus 0 curves. These may serve as a useful
introduction to this paper, since there are many similarities with
the integrals we study here, even though the geometry underlying
higher loop Feynman integrals is considerably more complex.

\subsection{Acknowledgements}
I owe many thanks to John Gracey and David Broadhurst for explaining
to me the problem of studying the three-loop contributions to the
massless two-point function, and also to Dirk Kreimer for many
helpful discussions.
The list of primitive divergent graphs up to 7 loops was kindly
given to me  by Karen Yeats, and my interest in the subject owes
much to the work of Bloch, Esnault and Kreimer \cite{B-E-K}, and the
encouragement of Pierre Cartier.

 This work was motivated by  the
fruitful workshop on perturbative approaches to quantum field theory
at ESI Vienna, and was continued at the MPIM during an IPDE
fellowship, and completed at the IHES.

\newpage

\section{Parametric representations}
Let $G$ be a Feynman graph with $h$ loops, $L$ internal edges,  and
$E$ external legs. The {\it graph polynomial} of $G$ is a homogeneous
polynomial of degree $h$ in variables $\alpha_1,\ldots, \alpha_L$
indexed by the set of internal edges of $G$. It is defined by the formula:
\begin{equation}\label{Udefn}
U_G = \sum_{T} \prod_{\ell\notin T} \alpha_\ell\ .
\end{equation}
The sum is over all spanning trees $T$ of $G$,  i.e., subgraphs
$T$ of $G$ which pass through every vertex of $G$ but which contain
no loop. Next consider the  homogeneous
polynomial of degree $h+1$
defined by:
\begin{equation}\label{Vdefn}
\V_G =  \sum_{S} \prod_{\ell\notin S} \alpha_\ell\, (q^{S})^2\ .
\end{equation}
The sum is over graphs $S\subset G$ with exactly two connected
components  $S=T_1\cup T_2$ where both $T_1$ and $T_2$ are trees,
such that $S$ is obtained by cutting a spanning tree $S'$ along an
edge $e$, and $q^S$ is the momentum flowing through $e$ in $S'$.
\begin{example}
 Let $G$ denote the two-point two-loop diagram
depicted in Fig $\ref{figureThreeLoop}$. Then
\begin{eqnarray} \label{U2p2L}
\qquad U_{G}\!\!\!&=\!\!\!&(\alpha_1+\alpha_5)(\alpha_2+\alpha_4) +
\alpha_3(\alpha_1+\alpha_2+\alpha_4+\alpha_5) \ ,  \\
\qquad
\V_{G}\!\!\!&=\!\!\!&\big(\alpha_3(\alpha_1+\alpha_2)(\alpha_4+\alpha_5)+
(\alpha_2\alpha_4\alpha_5+\alpha_1\alpha_4\alpha_5+\alpha_1\alpha_2\alpha_5+\alpha_1\alpha_2\alpha_4)\big)
\,q^2\ . \nonumber
\end{eqnarray}
\begin{figure}[h!]
  \begin{center}
    \epsfxsize=11.1cm \epsfbox{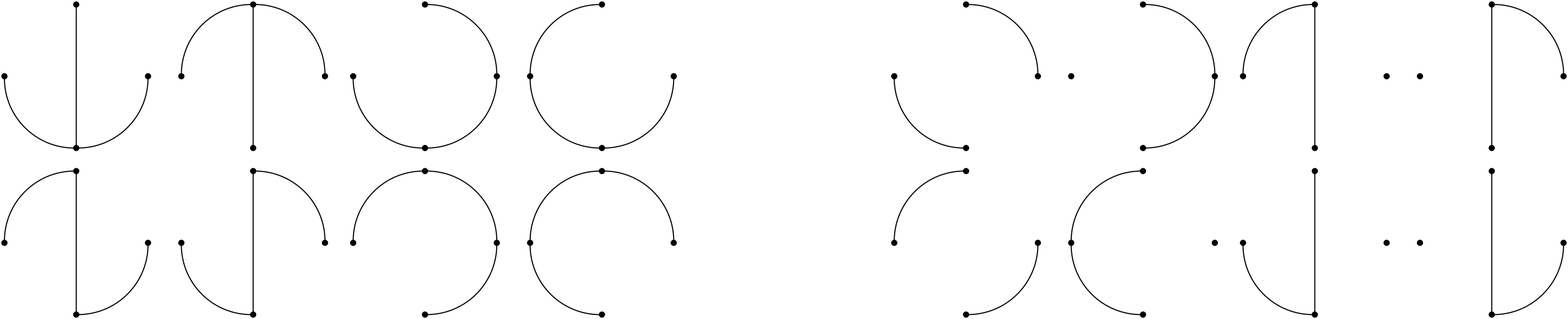}
  \label{figure1}
 \caption{On the left are shown the eight spanning trees for $G$, corresponding to the eight terms in $U_G$. On the right are the eight pairs of
trees $T_1\cup T_2$ which correspond to the eight terms in $\V_G$.}
  \end{center}
\end{figure}
\end{example}
To each internal edge $\ell$ of $G$ we associate a decoration
$a_\ell$, which is a positive real number (the power to which the
corresponding propagator is raised) and set
\begin{equation}\label{adef}
a=\sum_{1\leq \ell \leq L} a_\ell\ .\end{equation}
 The unregularised massless Feynman integral
$ I_G(\underline{a}, \underline{q},D)$ in dimension $D$, expressed
in Schwinger coordinates (\cite{IZ};\cite{Sm}, (2.36)) , is :
\begin{equation}\label{Schwinger}
(-1)^a\,{e^{i\pi[ a +h(1-D/2)]/2} \pi^{hD/2} \over \prod_{\ell}
\Gamma(a_\ell)}\int_0^\infty d\alpha_1 \ldots \int_{0}^\infty
d\alpha_\ell \,{\prod_{\ell} \alpha_\ell^{a_\ell-1} e^{i\V_G/U_G}
\over U_G^{D/2}} \ .
\end{equation}
 Now let $\lambda$ denote a non-empty set of internal edges
of $G$. By making the change of variables $\alpha_i=\alpha'_i\, t$,
for $1\leq i\leq L$,  where  $t=\sum_{\ell \in \lambda}
\alpha_\ell$, integrating with respect to $t$ from $0$ to $\infty$,
and finally replacing $\alpha'_i$ with $\alpha_i$ once again, this
integral can be  rewritten (\cite{Sm}, (3.32)):
\begin{equation}\label{Feynman}
{ \pi^{hD/2}\Gamma(a-hD/2) \over \prod_{\ell} \Gamma(a_\ell)}
\int_{H_\lambda}  \,{\prod_{\ell} \alpha_\ell^{a_\ell-1} \over
U_G^{D/2}} \Big({ U_G\over \V_G}\Big)^{a-hD/2} \, \Omega_L\ .
\end{equation}
Here, $\Omega_L=\sum_{i=1}^L (-1)^i\alpha_i d\alpha_1\ldots
d\widehat{\alpha_i} \ldots d \alpha_L$, and $H_\lambda=\{\alpha_i:
\sum_{l\in \lambda} \alpha_l = 1\}$. The fact that the integral does
not depend on the choice of subset $\lambda$ follows from the fact that it is really a projective integral, but is  known in the physics
literature as the  Cheng-Wu theorem. In this paper,  we shall always
take $\lambda$ to be a single edge, and set $D=4-2\,\varepsilon$.

\subsection{Primitive divergent graphs}
We say that  $G$ is \emph{broken primitive divergent} if it has
exactly 2 external legs, contains no divergent subgraphs, and
satisfies:
\begin{equation} \label{L2h}
L=2h+1\ .
\end{equation}
In this case, let $q$ denote the momentum entering or leaving each
external leg. By $(\ref{Vdefn})$, we  have $\V_G= V_G \, q^2$, where
$V_G$ is a polynomial in the $\alpha_i$ and does not depend on $q$.
Up to gamma-factors, we can therefore rewrite $(\ref{Feynman})$ as
the integral:
\begin{equation}\label{Idefn}
 I (a_1,\ldots, a_L, \varepsilon)= q^{h(2-\varepsilon)-a} \int_{H_\lambda} {\Omega_L \over U_G V_G} \prod_{\ell=1}^L \Big({\alpha_\ell U_G \over V_G}\Big)^{a_\ell-1} \Big({U_G^{h+1} \over V_G^h} \Big)^\varepsilon\ ,
\end{equation}
which converges for all $\Real a_\ell>0$. As is customary,  we set
\begin{equation} \label{nidef} a_\ell = 1 + n_\ell\, \varepsilon \ , \quad
\hbox{ for } 1\leq \ell \leq L\ ,
\end{equation}
where $n_\ell$ are positive integers. Before taking the Taylor
expansion with respect to $\varepsilon$, we first express the
integrand of $(\ref{Idefn})$ 
as an
integral of a simpler function.

When $G$ is broken primitive divergent, we can close up the two
external legs of $G$ to form a graph $\widetilde{G}$ with $h+1$
loops, and $L+1=2(h+1)$ edges (see figure 1, right). In
\cite{B-E-K}, such graphs were called primitive divergent,
and one
verifies that
\begin{equation} \label{UGtilde}
U_{\widetilde{G}}= V_G + \alpha_{L+1} U_G\ , \end{equation} where
$\alpha_{L+1}$ is the parameter attached to the edge of
$\widetilde{G}$ obtained by gluing the external legs of $G$ together
(see $\S\ref{sectnotation}$ below). The following lemma is nothing
other than the definition of Euler's beta function.

\begin{lem}Let $0< r< s$, and $u,v \neq 0$. Then
$$\int_{0}^{\infty} {x^{r-1} \over (u\,x+v)^s} dx = {1 \over u^{r}v^{s-r}}{ \Gamma(r)\Gamma(s-r) \over \Gamma(s)}\ .$$
\end{lem}
On comparison with $(\ref{Feynman})$, we  set $r=(h+1){D \over 2} -
a$, and $s={D\over 2}$. Substituting $(\ref{L2h})$ and
$(\ref{nidef})$ into the expression for $r$ leads us to define:
\begin{eqnarray}
a_{L+1} =  1 - \big( h+1 + \sum_{i=1}^L n_i \big) \, \varepsilon\ .
\end{eqnarray}
Therefore, using the previous lemma,  we can rewrite the integral
$(\ref{Idefn})$  as the product of certain explicit gamma factors
with the following integral:
\begin{eqnarray} \label{I-1defn}
 I (a_1,\ldots, a_{L+1}, \varepsilon)=  \int_{H_\lambda} {\prod_{\ell=1}^{L+1} \alpha_\ell^{a_\ell-1} \over  U_{\widetilde{G}}^{2-\varepsilon}}\,  \Omega_{L+1} \ ,
\end{eqnarray}
where $\Omega_{L+1} = \Omega_L d\alpha_{L+1}$,
and $H_\lambda$ is the hyperplane $\alpha_\lambda=1$ in
$[0,\infty]^{L+1}$.
%
%
The  Taylor expansion of the original Feynman integral
$(\ref{Feynman})$ can be retrieved from the following well-known
formula for the Taylor expansion of the gamma function:
$$\Gamma(s+1) =\exp( -\gamma s)\, \exp\Big( \sum_{k=2}^\infty (-1)^k\zeta(k) s^k\Big)\ . $$
The coefficients in the Taylor expansion of $(\ref{I-1defn})$ are
given by:
\begin{equation} \label{genint}
\int_{H_\lambda} {\Omega_{L+1} \over U^2_{\widetilde{G}}}
P(\alpha_i,\log(\alpha_i), \log(U_{\widetilde{G}}))\ ,
\end{equation}
where $m$ is an integer, and $P$ is a polynomial in $\alpha_i,
\log(\alpha_i)$ for $1\leq i\leq L+1$ and  $\log U_{\widetilde{G}}$
with rational coefficients. This is a period integral in the sense
of \cite{Ko-Za}.

In the first half of this paper, we will state a criterion on
$\widetilde{G}$ for $(\ref{genint})$ to evaluate to multiple zeta
values, and in the second half, we will outline an algorithm for
computing these values.

\subsection{Contraction and deletion of edges}\label{sectnotation}
 For any graph $G$,  its graph polynomial $U_G$ is linear with respect to each variable $\alpha_i$. We can write
\begin{equation}\label{Uedge}
U_G= U_{G\q \{i\}}+U_{G\backslash\{i\}}  \alpha_i\ ,\end{equation}
where $G\q\{i\}$ is the graph obtained from $G$ by contracting the
edge labelled $i$, and $G\backslash \{i\}$ is the graph obtained
from $G$ by deleting the edge labelled $i$.
 Using  $(\ref{UGtilde})$,
this implies that
\begin{equation} \label{UVsquish}
U_G=U_{\widetilde{G}\backslash\{L+1\}} \quad \hbox{ and } \quad V_G=U_{\widetilde{G}\q\{L+1\}}\ .
\end{equation}
From $(\ref{Uedge})$ and $(\ref{UVsquish})$, one deduces a similar
formula for $V_G$, namely:
\begin{equation}\label{Vedge}
V_G= V_{G\q \{i\}}+V_{G\backslash\{i\}}  \alpha_i\ .
\end{equation}

\begin{notation}
When the graph $G$ is implicit,  we will often write $U,V$ instead
of $U_G, V_G$. In order to lighten the notation, we   adopt the
useful convention from \cite{B-E-K}, which consists in writing,  for
any polynomial $\Psi$ which is linear in the variable $\alpha_i$:
$$\Psi^{(i)} = {\partial \over \partial \alpha_i} \Psi \ , \quad \hbox{ and }  \quad \Psi_i = \Psi\Big|_{\alpha_i=0}\ .$$
We   therefore have  $U_{i}=U_{G\q \{i\}}$,
$U^{(i)}=U_{G\backslash\{i\}}$, and $(\ref{Uedge})$ and
$(\ref{Vedge})$ become
$$U=U_i + U^{(i)} \alpha_i \quad \hbox{ and } \quad  V=V_i + V^{(i)} \alpha_i \ .$$
Since the operations of contracting and deleting distinct edges
commute, we can write $U_{12} = U_{G\q\{1,2\}}$,
$U^{(1)}_2=U_{G\q \{2\}\backslash \{1\}}$, $U^{(12)}= U_{G\backslash
\{1,2\}}$, and so on, where indices in the superscript (subscript)
correspond to deleted (resp. contracted) edges. \end{notation}

\newpage
\section{Symbolic integration using polylogarithms}
\label{sectSymbolicInt}
 Let $G$ be a broken primitive divergent graph, and let $\widetilde{G}$ be the graph obtained by closing its external legs.
To illustrate our integration method, let us begin to compute the
integral $(\ref{Idefn})$ in the case where all decorations
$a_\ell=1$, and $\varepsilon=0$. Therefore consider  the
convergent integral
\begin{equation}\label{IGminus1}
I_{\widetilde{G}} = \int_{H_\lambda} {\Omega_{L+1} \over
U_{\widetilde{G}}^2}\ , \qquad \qquad \qquad \hbox{{\it (Step 0)}}
\end{equation}  
which was studied in \cite{B-E-K}.  We will assume that $\lambda$ is
a single index, say $L$. In this case we can write $\Omega_{L} =
d\alpha_1 d\alpha_2\ldots d\alpha_{L-1}$ and $\Omega_{L+1}= \Omega_L
d\alpha_{L+1}$. The domain of integration is simply
$H_\lambda=\{\alpha_{L}=1,\, 0\leq \alpha_i \leq \infty\ , i\neq
L\}$. Using equation $(\ref{UGtilde})$ to replace
$U_{\widetilde{G}}$ with $V_G+ \alpha_{L+1}U_G$, we can perform one
integration with respect to  $\alpha_{L+1}$ from $0$ to $\infty$, to
obtain:
\begin{equation}\label{IG0}
I_G= \int_{\alpha_L=1} \int_0^\infty {d\alpha_{L+1} \over (V_G+\alpha_{L+1}U_G)^2} \,\Omega_L= \int_{H_\lambda} {\Omega_L \over U_G V_G}  \qquad \qquad \qquad \hbox{{\it (Step 1)}} 
 \ .\end{equation}
From here on, we write $U, V$ instead of $U_G, V_G$ and use the
notations of $\S\ref{sectnotation}$.  Now choose any variable in the
integrand, say $\alpha_1$, and write $I_G$ in terms of $\alpha_1$:
$$I_G= \int_{\alpha_L=1} {\Omega_L \over (U_1 +U^{(1)} \alpha_1 )( V_1 +V^{(1)} \alpha_1 )}\ .$$
By decomposing into partial fractions with respect to $\alpha_1$,
this gives
$$I_G= \int_{H_\lambda} {\Omega_L \over U^{(1)}V_{1} -  U_1 V^{(1)}  }
 \Big({U^{(1)} \over  U_1 + U^{(1)} \alpha_1 }-{V^{(1)} \over  V_1 + V^{(1)}\alpha_1 }\Big)\ .$$
Now we can peform an integration with respect to $\alpha_1$ from $0$ to $\infty$ at the expense of introducing a new function, namely the logarithm.
\begin{equation}\label{IG1}
I_G= \int_{\alpha_L=1} {\log U^{(1)} -\log U_{1} -\log V^{(1)} +\log
V_{1} \over U^{(1)}V_1 - U_1 V^{(1)} } \, \,d\alpha_2\ldots
d\alpha_{L-1}  \qquad \hbox{{\it (Step 2)}} \ .
\end{equation} The logarithm  $\log f$ should be
regarded as a symbol which satisfies the formal rule $d\log f =
f^{-1} df$, since changing the constant of integration does not
affect  the integrand of $(\ref{IG1})$. At this stage, something
remarkable occurs.  The denominator factorizes as the square of a
polynomial $D$ which is linear in each variable
 $\alpha_2,\ldots, \alpha_L$:
 \begin{equation} \label{Ddefn}
D^2=U^{(1)}V_{1} -  U_{1} V^{(1)}
\end{equation}
As observed in \cite{B-E-K}, this phenomenon is quite general and
follows from a result due to Dodgson on determinants of matrices.
One can also give a formula for $D$ in terms of trees in $G$, but
this will not be required here.

We can repeat this argument. Choose a variable, say $\alpha_2$, and write
$$D=D_2\alpha_2+ D^{(2)}\ .$$
By decomposing the integrand of $(\ref{IG1})$ into partial fractions
with respect to the variable $\alpha_2$, we can integrate with
respect to the variable $\alpha_2$. To simplify the  notation, we
define
\begin{equation}\label{defCurlynotation}
\{p,q|r,s\}={p\big(\log(p)+\log(s)-\log(q)-\log(r)\big) \over r(ps-qr) }+ {\log(q) \over rs}\ .\end{equation}
One verifies that
$$\int_{0}^\infty {\log(px+ q) \over (rx+s)^2} dx = \{p,q|r,s\}\ .$$

\begin{cor} \label{corstep3} The integral $I_G$ is equal to the $L-2$ dimensional integral:
$$I_G= \int_{\alpha_L=1} \Big(\{ U^{(1)}_{2}, U^{(1,2)} | D_2, D^{(2)} \}-\{ U_{12}, U^{(2)}_1 | D_2, D^{(2)} \}  - \{ V^{(1)}_{2}, V^{(1,2)} | D_2, D^{(2)}  \}$$
$$+\{ V_{12}, V^{(2)}_1 | D_2, D^{(2)} \}\Big) \,d\alpha_3\ldots d\alpha_{L-1}  \qquad \hbox{{\it (Step 3)}}\ . $$
\end{cor}

As long as there exists a variable with respect to which all polynomials which occur in the integrand are linear, this can be repeated.

 At the next stage of the integration
process, one has to introduce the dilogarithm\footnote{As we shall
see later, it is more convenient to use  the function satisfying the
differential equation $dL(f)=f^{-1}\log(f+1)\,df$, which is  a close
relative of the dilogarithm.}, which is formally defined by the
differential equation
\begin{equation} d \,\Li_2(f) = -{\log(1-f) \over f}\,df\
.\end{equation}

 In this manner, we obtain a conditional algorithm
for computing integrals such as $(\ref{genint})$, which can be
approximately formalised as follows:

\begin{enumerate}
  \item Choose a variable in which  all terms 
of the integrand are linear.
  \item Formally take a primitive of the integrand with respect to this variable.
  \item Evaluate this primitive  at 0 and $\infty$, and repeat.
\end{enumerate}
The algorithm fails exactly when in $(1)$, we can no longer find a variable with respect to which all terms in the integrand are linear.
This will be formulated precisely in the following section. It can be checked in advance whether the algorithm will terminate, as this
depends only on the topology of the graph $G$.
 When the algorithm does terminate,  any integral $(\ref{Idefn})$ with arbitrary decorations can
always be computed in a finite number of
steps in terms of a fixed differential algebra of polylogarithms determined in advance by the topology of $G$.
\\

For example, in the above, we can keep track of the singularities of
the integrand (i.e., terms occuring in the denominator, or arguments
of the logarithm functions) at each step. At the first 
step, we represent the  singularities as the set:
$$ \{U,V\}\ .$$
After the second integration, the singularities are given by:
$$ \{U^{(1)},U_1,V^{(1)},V_1,D\}\ .$$
At the third  step (corollary $\ref{corstep3}$), we have in the same
manner
$$ \{U^{(12)},U^{(1)}_2,U^{(2)}_1, U_{12} ,V^{(12)},V^{(1)}_2,V^{(2)}_1,V_{12},D^{(2)},D_2 \}\ ,$$
along with the set of irreducible factors of the four denominator
terms:
$$D^{(2)}U^{(1)}_2-D_2U^{(12)}\ , \ D^{(2)}U_{12}-D_2U^{(2)}_1 \ ,
\ D^{(2)}V^{(1)}_2-D_2V^{(12)} \ , \  D^{(2)}V_{12}-D_2V^{(2)}_1\ .$$ We are
therefore  led to consider an algorithm for the reduction of sets of
polynomials with respect to their variables. This will ultimately
lead  to a criterion for the computability of broken primitive
divergent Feynman graphs (theorem $\ref{thmMAINCRITERION}$).

\newpage
\section{ Reduction of polynomials}\label{sectFibAlg}
Keeping track of the polynomials, or singular loci, which can occur
in this integration process gives rise to  a  reduction
algorithm  which can be used to check the outcome of the integration
process without actually doing it.

\subsection{The simple reduction algorithm} 
 Let $S=\{f_1,\ldots, f_N\}$, where $f_1,\ldots, f_N$ are
polynomials in the variables $\alpha_1,\ldots, \alpha_m$, with
rational coefficients.
\begin{enumerate}
  \item Suppose that there exists an index $1\leq r\leq m$ with
respect to which every polynomial $f_1,\ldots, f_N$ is linear in the
variable $\alpha_r$. Then we can write:
$$f_i = g_i \, \alpha_r+ h_i \qquad \hbox{ for }\quad  1\leq i \leq N\ ,$$
where $g_i= \partial f_i/\partial \alpha_r$, and $h_i =
f_i|_{\alpha_r=0}$.
Define a new set of polynomials:
$$\widetilde{S}_{(r)}=\{ \big(g_i\big)_{1\leq i\leq N}\ , \  \big(h_i\big)_{1\leq i\leq N}\ ,\   \big(h_i g_j- g_i h_j\big)_{1\leq i<j\leq N}
\}\ .$$
  \item  Let $S_{(r)}$ be the set of irreducible factors of
polynomials in $\widetilde{S}_{(r)}$.
\end{enumerate}
The polynomials now occurring in $S_{(r)}$  are  functions of one
fewer variables, namely $\alpha_1,\ldots, \alpha_{r-1},\alpha_{r+1},
\ldots, \alpha_m$. This process can be repeated. If at each stage
there exists a variable in  which all polynomials are linear, we can
proceed to the next stage. This gives a sequence of variables
$\alpha_{r_1},\alpha_{r_2},\ldots,\alpha_{r_n}$, and a sequence of
sets
\begin{equation}\label{simpleseq}
 S_{(r_1)}\ , \ S_{(r_1,r_2)}\ ,\  \ldots\ ,\  S_{(r_1,\ldots,
r_n)}\ .
\end{equation}
 If there  exists a sequence $(r_1,\ldots, r_m)$ such that
every variable is eventually eliminated, then we say that the
reduction terminates. When this happens, we say that the set $S$ is
\emph{simply reducible}.

\begin{rem} \label{remCONVENTION}
We can remove any constants, and any monomials of the form
$\alpha_i$ which  occur as elements in $S_{(r)}$, as this does not affect the
outcome of the algorithm.
\end{rem}


\begin{defn}
Let $G$ be a broken primitive divergent graph, and let
$$S_G=\{U_{\widetilde{G}}\}\ .$$
 We say that
$G$ is simply reducible if $S_G$ is simply reducible.
\end{defn}
 Observe that $\{U_{\widetilde{G}}\}_{(L+1)}=\{U_G,V_G\}$, where $\alpha_{L+1}$ is the edge
variable of $\widetilde{G}$ obtained by closing  the external legs
of $G$, by equation   $(\ref{UGtilde})$. 


\begin{example} \label{example2loop}
Consider  the 2-loop 2-point graph $G$ depicted in fig. 1. We write
$S=S_G=\{U_{\widetilde{G}}\}$. After reducing with respect to
$\alpha_6$, we have $S_{(6)}=\{U,V\}$, where $U=U_G$, $V=V_G$ are
given by $(\ref{U2p2L})$:
\begin{eqnarray} \nonumber
S_{(6)}\!\!\!&=\!\!\!&\{
\alpha_1\alpha_2+\alpha_1\alpha_4+\alpha_5\alpha_2+\alpha_4\alpha_5
+
\alpha_3\alpha_1+\alpha_3\alpha_2+\alpha_3\alpha_4+\alpha_3\alpha_5\ ,  \\
&\quad
&\alpha_3\alpha_1\alpha_4+\alpha_3\alpha_1\alpha_5+\alpha_3\alpha_2\alpha_4+\alpha_3\alpha_2\alpha_5+
\alpha_2\alpha_4\alpha_5+\alpha_1\alpha_4\alpha_5+\alpha_1\alpha_2\alpha_5+\alpha_1\alpha_2\alpha_4
 \}\ .\nonumber
\end{eqnarray}
Since both polynomials are linear in $\alpha_1$, we can reduce with
respect to $\alpha_1$ to obtain
$$\widetilde{S}_{(6,1)}=\{U^{(1)},U_1,
V^{(1)},V_1, U^{(1)}V_1-U_1V^{(1)}\}\ .$$ But
$U^{(1)}V_1-U_1V^{(1)}$ factorizes, by the Dodgson identity, and so
at the second stage we have  $S_{(6,1)}=\{U^{(1)},U_1, V^{(1)},V_1,
D\}$, where $D^2=(U^{(1)}V_1-U_1V^{(1)})$: 
\begin{eqnarray} \nonumber
 S_{(6,1)}=\{& \!\!\!\!\!\!\alpha_2\alpha_5 + \alpha_4\alpha_5+\alpha_2\alpha_3+\alpha_3\alpha_4+\alpha_3\alpha_5\ ,
\,\,\alpha_2+\alpha_3+\alpha_4\ ,\,\,
  \alpha_3\alpha_4+\alpha_3\alpha_5+\alpha_4\alpha_5 \ ,\,\, &\!\!\!\!   \\
 & \!\!\!\!\nonumber\alpha_3\alpha_4+\alpha_3\alpha_5+\alpha_4\alpha_5+\alpha_2\alpha_5+\alpha_2\alpha_4\ , \,\, \nonumber\alpha_3\alpha_4+\alpha_4\alpha_5+\alpha_2\alpha_5+\alpha_3\alpha_5\, \} &
\end{eqnarray}
Since each polynomial is linear, we can reduce with respect to the
variable $\alpha_2$:
$$
 S_{(6,1,2)}=\{ \alpha_4+\alpha_5\ ,\,\, \alpha_3\alpha_4 + \alpha_3\alpha_5+\alpha_4\alpha_5\
,\,\, \alpha_3+\alpha_4\ ,\,\,
  \alpha_3+\alpha_5 \ ,\,\, \alpha_3-\alpha_4  \}\ .
$$
Next, reducing with respect to the variable $\alpha_5$, gives
$$
 S_{(6,1,2,5)}=\{ \alpha_3+\alpha_4\ , \,\,\alpha_3-\alpha_4\}\ ,
$$
and finally, reducing with respect to $\alpha_3$ gives
$$
S_{(6,1,2,5,3)}=\emptyset\ .
$$
Here, and from now on, we will adopt the convention (remark
$\ref{remCONVENTION}$) that we remove all constant terms and
variables $\alpha_i$ from each stage of the reduction.
\end{example}

In $\S\ref{sectSymbolicInt}$, the domain of integration for $I_G$
was  a hyperplane $\alpha_\lambda=1$ for some index $\lambda$. If
$S_G$ is simply reducible with respect to some order
$(\alpha_{r_1},\ldots, \alpha_{r_n})$ of the variables, we can
choose $\lambda$ to be the index of one of the two final variables
$\alpha_{r_{n-1}},\alpha_{r_n}$. In the previous example, we can
take $\lambda=4$, and set $\alpha_4=1$, $S'=S_G\big|_{\alpha_4=1}$.
This gives
$$
 S'_{(6,1,2)}=\{ \alpha_5+1\ ,\,\, \alpha_3+\alpha_5 + \alpha_3\alpha_5\
,\,\, \alpha_3+1\ ,\,\,
  \alpha_3+\alpha_5 \ ,\,\, \alpha_3-1  \}\ ,
$$
and
$$ S'_{(6,1,2,5)}=\{ \alpha_3+1\ , \,\,\alpha_3-1\}\ .$$
The reduction algorithm reflects the set of singularities which
occur at each stage of the integration process. Roughly speaking,
the integrand at the $k^\mathrm{th}$ stage, after integrating with
respect to $\alpha_{r_1},\ldots, \alpha_{r_k}$, will have
singularities along the zero locus of polynomials in the sets
$S_{(r_1,\ldots,r_m)}$, for $m\geq k$, and along  the  axes
$\alpha_i=0$.

\begin{rem} The  algorithm described above is  similar to an algorithm defined by
Stembridge to study the zeros of graph polynomials over finite
fields. A similar argument was  also used in  \cite{B-E-K} to obtain
a Tate filtration on  graph hypersurfaces.
\end{rem}

In the previous example, we see that at the penultimate stage, we
expect to have polylogarithms in $\alpha_3$, with singularities
along $\alpha_3=0, \alpha_3=-1,\alpha_3=1$. This  would imply that
at the final stage, the  integral $I_G$ is a rational linear
combination of alternating multiple sums $\MZV^2$, i.e., periods of
$\Pro^1\backslash \{0,-1,1,\infty\}$. This is not good enough, since
it is known that the  integral $I_G$ in fact gives multiple zeta
values. To rectify this problem, we introduce the Fubini reduction
algorithm.


\subsection{The Fubini reduction algorithm} \label{sectStrongRed}
In the simple reduction algorithm, the sets $S_{(r_1,\ldots, r_k)}$
which control the singularities of the integration process,  depend
in an essential way on the order of the variables $r_1,\ldots, r_k$
which was chosen. However, it is an obvious consequence of Fubini's
theorem that the final integral does not depend on the particular
order of integration.  More precisely, if $F$ is the integrand
obtained at the $k^{\mathrm{th}}$ stage, we clearly have
$$\int_{0}^\infty d\alpha_{r_i} \int_{0}^\infty d\alpha_{r_{j}} F =
\int_{0}^\infty d\alpha_{r_j} \int_{0}^\infty d\alpha_{r_{i}}F\ . $$
The left-hand integral has singularities contained in the zero locus
of polynomials in  $S_{(r_1,\ldots, r_{k}, r_i, r_{j} )}$; the
right-hand integral has singularities contained in the zero locus of
polynomials in $S_{(r_1,\ldots, r_{k}, r_{j}, r_i)}$. It follows
that both have singularities contained in the zero locus of
polynomials in the intersection $S_{(r_1,\ldots, r_{k}, r_i, r_{j}
)}\cap S_{(r_1,\ldots, r_{k}, r_{j}, r_i)}. $

In general, we define sets recursively as follows:
\begin{eqnarray}\label{Ssquare}
S_{[r_1, r_2]} &= &S_{(r_1,r_2)} \cap S_{(r_2,r_1)} \ ,\\
S_{[r_1, r_2,\ldots, r_k]} &= &\bigcap_{1\leq i\leq k}
S_{[r_1,\ldots,\widehat{r_i}, \ldots, r_{k}](r_i)}\ , \quad k\geq 3\
, \nonumber
\end{eqnarray}
where in  the second line one applies a reduction to the set
$S_{[r_1,\ldots,\widehat{r_i}, \ldots, r_{k}]}$ with respect to the
variable $\alpha_{r_i}$ (steps (1) and (2) in the previous section).
It may happen that it is not possible to apply the reduction to
$S_{[r_1,\ldots,\widehat{r_i}, \ldots, r_{k}]}$   if it contains a
polynomial which is non-linear in $\alpha_{r_i}.$ In that case,
$S_{[r_1,\ldots,\widehat{r_i}, \ldots, r_{k}](r_i)}$ is undefined
and we omit it from the intersection in the second line of
$(\ref{Ssquare})$. It is also possible that this happens for all
$1\leq i\leq k$, in which case $S[r_1,\ldots, r_k]$ is undefined.

 This is the \emph{Fubini reduction algorithm},
and it gives rise to a  sequence of sets
\begin{equation}\label{Fubiniseq}
 S_{(r_1)}\ , \ S_{[r_1,r_2]}\ ,\  \ldots\ ,\  S_{[r_1,\ldots,
r_n]}\ .
\end{equation}

\begin{defn}
We say that $S$ is \emph{Fubini reducible} if  there  is a sequence
$(r_1,\ldots, r_m)$ such that every variable is eventually
eliminated, and such that every polynomial in $S_{[r_1,\ldots,
r_k]}$ is linear in $\alpha_{r_{k+1}}$, for all $1\leq k\leq m$.
\end{defn}
  If $G$  is a broken
primitive divergent graph,
 we say that
$G$ is Fubini reducible if $S_G=\{U_{\widetilde{G}}\}$ is Fubini
reducible.
 As in remark $\ref{remCONVENTION}$,
we
  can remove  constants and terms of the form
$\alpha_i$ without affecting  the algorithm.



\begin{rem}

By the definition $(\ref{Ssquare})$, the set $S_{[r_1, r_2,\ldots,
r_k]}$ is contained in the set $S_{[r_{\pi(1)}, r_{\pi(2)},\ldots,
r_{\pi(k)}]}$ for any permutation $\pi$ of $1,2,\ldots, k$. Thus if
$S$ is simply reducible, then the sequence $(\ref{Fubiniseq})$ is
contained in  the sequence $(\ref{simpleseq})$.
 In  general, the  inclusion
$$S_{[r_1, r_2,\ldots, r_k]} \subseteq \bigcap_{\pi \in \Sym_k} S_{(\pi(r_1),\ldots,
\pi(r_k))}\ ,$$ is strict, where the intersection is taken over all
permutations of $r_1,\ldots, r_k$. In practice, the cascade of
polynomials one obtains using the Fubini reduction  is considerably
smaller than with the simple reduction algorithm.
\end{rem}


\begin{example} \label{exFubinired}
We  compute the Fubini reduction algorithm for  the two-loop two
point function, as in example $\ref{example2loop}$. Let
$S=\{U_{\widetilde{G}}\}$, where $U_{\widetilde{G}}=\alpha_6\, U_G+
V_G$, and $U_G,V_G$ are given by $(\ref{U2p2L})$.
 One checks  that $ S_{[2,6](5)}$ is given by:
$$ \{ \alpha_3+\alpha_4\  , \alpha_3+\alpha_4 +\alpha_1\ ,
\alpha_1\alpha_4+\alpha_1\alpha_3+\alpha_3\alpha_4\ ,
\alpha_1+\alpha_3\ ,
  2\alpha_3\alpha_4+\alpha_3^2+\alpha_1\alpha_4+\alpha_1\alpha_3 \ , \alpha_1\alpha_3+\alpha_1\alpha_4+\alpha_3^2 \}
$$
and $
 S_{[5,6](2)}$ by
$$ \{ \alpha_3+\alpha_4\  , \alpha_3+\alpha_4 +\alpha_1\ ,
\alpha_1\alpha_4+\alpha_1\alpha_3+\alpha_3\alpha_4\ ,
\alpha_1+\alpha_3\ ,
  \alpha_1\alpha_4+\alpha_3^2+\alpha_3\alpha_4 \ , 2 \alpha_1\alpha_3+\alpha_1\alpha_4+\alpha_3^2 +\alpha_3\alpha_4\}
$$
and $
 S_{[2,5](6)}$ by
$$ \{ \alpha_3+\alpha_4\  , \alpha_3+\alpha_4 +\alpha_1\ ,
\alpha_1\alpha_4+\alpha_1\alpha_3+\alpha_3\alpha_4\ ,
\alpha_1+\alpha_3 \ , \alpha_1+\alpha_4, \alpha_1-\alpha_4\}\ .
$$
Taking the intersection of all three sets gives
$$
 S_{[2,5,6]}=\{ \alpha_3+\alpha_4\  , \ \alpha_3+\alpha_4 +\alpha_1\ ,\ \alpha_1\alpha_4+\alpha_1\alpha_3+\alpha_3\alpha_4\ ,\
\alpha_1+\alpha_3\}\ .
$$
By performing an ordinary reduction with respect to the variable
$\alpha_1$, one obtains:
$$S_{[2,5,6](1)}=\{\alpha_3+\alpha_4, \alpha_3
\alpha_4+\alpha_3^2+\alpha_4^2\}\ .$$
By computing $S_{[6,1,2](5)}$, $S_{[6,1,5](2)}$, $S_{[1,2,5](6)}$
and intersecting all four (actually it suffices to intersect with
$S_{(6,1,2,5)}$ from example $\ref{example2loop}$ in this case), one
verifies that
$$S_{[6,1,2,5]}= \{\alpha_3+\alpha_4\}\ .$$
\end{example}

As before,  we can take our domain of integration to be the
hyperplane $\alpha_4=1$. At the final stage of the integration
process, the Fubini reduction algorithm will predict (see theorem
$\ref{thmMAINCRITERION}$) that the integrand has singularities in
$\alpha_3=0, \alpha_3=1$. Therefore we expect to obtain multiple
polylogarithms with singularities in $\{0,1,\infty\}$ at the
penultimate stage, and hence multiple zeta values as the final
answer. In this way, the problem raised at the end of the previous
section has apparently been overcome.
\\

However, constants can appear at every stage during the integration
process (see \S7), and one needs to verify that we only obtain
multiple zeta values every time. This leads to a further
ramification condition to be verified, for each set $S_{[r_1,\ldots,
r_k]}$.


\subsection{The ramification condition}\label{sectRamcond}
Let $S$ be Fubini reducible for some order $(r_1,\ldots, r_m)$ of
the variables. After setting the final variable $\alpha_{r_m}=1$  in
the sequence  $(\ref{Fubiniseq})$, we obtain  a new sequence of sets
(the reductions of $S'=S\big|_{\alpha_{r_m}=1}$):
\begin{equation}\label{Fubiniseqred}
 S'_{(r_1)}\ , \ S'_{[r_1,r_2]}\ ,\ \ldots\ ,\
S'_{[r_1,\ldots, r_{m-1}]}\ ,
\end{equation}
 where every polynomial in $S'_{[r_1,\ldots,
r_k]}$ is linear in the variable $\alpha_{r_{k+1}}$.  Therefore, if
we write  $ S'_{[r_1,\ldots, r_k]}=\{f_1,\ldots, f_{M_k}\}$, then we
have
$$ f_i= a_i\, \alpha_{r_{k+1}} + b_i\ , \quad \hbox{for } 1\leq i\leq M_k\ ,$$
where $a_i, b_i$ are polynomials in $\alpha_{r_{k+2}},\ldots,
\alpha_{r_{m}}$. We define  $\Sigma_{\alpha_k}$ to be
\begin{equation}\label{firstsigmakdef}
\Sigma_{\alpha_k}= \{ -{b_i \over a_i}\, \hbox{ such that } a_i\neq
0 \}\ .\end{equation} The set $\Sigma_{\alpha_k}$ clearly depends on
the ordering of the variables $(r_1,\ldots, r_m)$.
\begin{defn}\label{defnunram}
We say that $\Sigma_{\alpha_k}$ is unramified if:
\begin{equation} \label{ramcond}
\lim_{\alpha_{r_m}\rightarrow 0}
\Big(\lim_{\alpha_{r_{m-1}}\rightarrow 0}\Big(\ldots
\Big(\lim_{\alpha_{r_{k+2}}\rightarrow
0}\Sigma_{\alpha_k}\Big)\Big)\cdots\Big) \subseteq \{0,-1,\infty\} \ .
\end{equation}
We say that the sequence of sets $(\ref{Fubiniseq})$  is unramified
if the corresponding $\Sigma_{\alpha_k}$ are unramified for all
$1\leq k\leq m-1.$ Finally, we say that $S$ is unramified if it is
Fubini reducible, and if there exists a sequence $(r_1,\ldots,
r_m)$ such that $(\ref{Fubiniseqred})$ is unramified.
\end{defn}

\begin{example} One can check that the sequence of sets given in example
$\ref{exFubinired}$ is unramified. We check the last two terms only.
Setting $\alpha_4=1$, we have:
$$
 S'_{[6,2,5]}= \{ 1+\alpha_3\  , \ 1+\alpha_1 +\alpha_3\ ,\ \alpha_1+\alpha_3+\alpha_1\alpha_3\ ,\
\alpha_1+\alpha_3\}\ ,\ S'_{[6,1,2,5]}= \{1+\alpha_3\}\ .$$ Thus
$$\Sigma_{\alpha_1} = \{ -(\alpha_3+1)\ , \  -{\alpha_3\over
1+\alpha_3}\ ,\  -\alpha_3\}\ , \hbox{ and }
\Sigma_{\alpha_3}=\{-1\}\ .$$ The ramification condition
$(\ref{ramcond})$ is  satisfied in  both cases, since
 $\Sigma_{\alpha_3}\subset \{0,-1,\infty\}$, and
$\lim_{\alpha_3\rightarrow 0} \Sigma_{\alpha_1} = \{0,-1\}$.
\end{example}

\subsection{The main theorem} \label{sectMainthm}
The main theorem gives a simple criterion for a master Feynman
integral to evaluate to multiple zeta values.

\begin{thm}\label{thmMAINCRITERION} Let $G$ be a broken primitive divergent Feynman diagram
and let $S_{G}=\{U_{\widetilde{G}}\}$. If $S_{G}$ is Fubini
reducible and unramified, then the coefficients in the Taylor
expansion of the integral $I(G)$ $(\ref{I-1defn})$ are rational
linear combinations of multiple zeta values.
\end{thm}

We also require a variant to allow for multiple zeta values ramified
at roots of unity. Let $S$ be Fubini reducible, as above. We say
that $S$ is  ramified at $p^{\mathrm{th}}$ roots of unity if  there
exists a sequence $(r_1,\ldots, r_m)$ such that the corresponding
sets $\Sigma_{\alpha_k}$ satisfy    
$$\Sigma_{\alpha_{k}} \subseteq \{0,\infty\}\cup \{ -\omega: \omega^p=1\} \ .
$$
If $p=1$, then this coincides with the definition $\ref{defnunram}$
above.

\begin{thm}\label{thmMAINCRITERION2} Let $G$ be a broken primitive divergent Feynman diagram
and let $S_{G}=\{U_{\widetilde{G}}\}$. If $S_{G}$ is Fubini
reducible and ramified at $p^{\mathrm{th}}$ roots of unity, then the
coefficients in the Taylor expansion of the integral $I(G)$ are
rational linear combinations of multiple zeta values ramified at
$p^{\mathrm{th}}$ roots of unity, {\it i.e.}, in $\MZV^p$.
\end{thm}

The results stated in the introduction are the result of applying
these two theorems to the set of all primitive divergent graphs up
to 5 loops. We computed the Fubini reduction algorithm for the
cross-hairs diagram (for which the graph polynomial has 45 terms),
the graphs $5\,R$, $5\,P$ and $5\,N$ (for which it has 128,130
and 135 terms, respectively), and all eight primitive divergent graphs with 6 loops and crossing number 0 and 1, as depicted in figures
6 and 7 (for which the graph polynomials have around 400 terms).
We found  that all planar bpd graphs up to 5
loops are Fubini reducible and unramified, and all non-planar bpd graphs with crossing number 1 are Fubini reducible but ramified
 at $6^\mathrm{th}$
roots of unity. 
The bipartite graph $K_{3,4}$ has crossing number 2 (figure 8), and  is not Fubini reducible.
These calculations were done using maple.

The outcome of the Fubini reduction algorithm may depend on the particular choice of two final variables. Typically, the ramification
is better if one chooses the two final variables to correspond to edges meeting at a four-valent vertex.
In each example, we  computed the algorithm for all possible pairs of  final variables, and chose a pair giving the least ramification. In an ideal world,
one could simply take the intersection of the period rings obtained for every such choice, and the algorithm would in that case not depend on any choices.
In general, however, taking the intersection of rings of periods  requires some powerful diophantine results, which are at present totally out of reach.

\subsection{Extensions} \label{sectExtensions}
There are a number of ways in which the reduction algorithm might be improved.
First of all,  in our computer calculations, we factorized our polynomials over $\Q$, rather than over the algebraic
closure $\overline{\Q}$, which may or may not have made a significant difference.

More interestingly, the Fubini reduction method stops as soon as it finds that, for every variable $\alpha_i$,
 there is a polynomial which is quadratic in that variable. It is possible, by introducing new variables which are the square roots of discriminants (by passing
to a ramified cover), one could generalize  the method
 to deal with some of the quadratic terms. The basic idea is that, in the case of a plane curve for example, a polynomial need not necessarily be
linear to define a curve of genus 0. A conic, of degree 2,  also defines a curve of genus 0, and our method should  also extend to this case.
A genuine obstruction should occur at the degree 3 level, since one expects to find elliptic components in this case. In any case, by extending our method
 to deal with quadratic terms, a larger class of graphs may then become reducible, and amenable to computation.

Finally, if one is interested only in the leading term of the Taylor expansion of the graph $I(G)$, then one can simplify the reduction algorithm. Recall that,
for the general
integral (with arbitrary logarithms in the numerator), we obtained at the second step, the set of singularities:
$$S_{[L+1,1]} = \{U^{(1)},U_1,V^{(1)},V_1,D\}\ .$$
Notice  that in $(\ref{IG1})$, where  the numerator has no logarithmic terms,    not all pairs of singularities can occur. Thus one should perform a
Fubini reduction separately on the four sets
$$\{U^{(1)},D\} \ , \ \{U_1,D\} \ , \ \{V^{(1)},D\} \ , \ \{V_1,D\} \ ,$$
which may improve the ramification. There may  also exist graphs which are not Fubini reducible, but for which each of the four
sets above is Fubini reducible. This would prove  that the  first term only in the Taylor expansion is a multiple zeta value.
\\

In the remainder of the paper we explain how to compute the
coefficients in the Taylor expansion when the conditions of  theorems $\ref{thmMAINCRITERION}$ and $\ref{thmMAINCRITERION2}$ hold.

\newpage
\section{Hyperlogarithms, polylogarithms, and primitives}

In this section, we outline the function theory underlying the
iterated integration procedure. For a more detailed introduction,
see \cite{Br1, Br2, Br3}.
\subsection{Hyperlogarithms} \label{sectHyperlogarithms}
Let $\Sigma=\{\sigma_0,\sigma_1,\ldots, \sigma_N\}$, where
$\sigma_i$ are distinct points of $\C$. We will always assume that
$\sigma_0=0$. For now we will consider the points $\sigma_i$ to be
stationary, but later we will allow them to be variables in their
own right.

Let $A=\{\ao_0,\ao_1,\ldots, \ao_N\}$ denote an alphabet on $N+1$
letters, where each symbol $\ao_i$ corresponds to the point
$\sigma_i$. Let $\Q\langle A\rangle$ denote the vector space
generated by all words $w$ in the alphabet $A$, along with the empty
word which we denote by $e$. To each such word $w$, we associate a
hyperlogarithm function:
$$L_w(z) : \C \backslash \Sigma \To \C$$
which is  multivalued, i.e., it is  a meromorphic function on the
universal covering space of $\C \backslash \Sigma$. Let $\log(z)$
denote the principal branch of the logarithm.

\begin{defn}
\label{defnHyperlog} Let $A^\times$ denote the set of all words $w$
in $A$, including $e$. The family of functions $L_w(z)$ is uniquely
determined by the following properties:
\begin{enumerate}
  \item $L_e(z) = 1$, and  $L_{\ao_0^n}(z) = {1\over n!} \log^n(z)$ , for all $n\geq 1$.
  \item For all words $w\in A^\times$, and all $0\leq i\leq N$,
 $${\partial \over \partial z} L_{\ao_i w} (z) = {1 \over
z-\sigma_i} L_w(z) \ ,\  \hbox{ for } z\in \C\backslash \Sigma\ .$$
  \item For all words $w\in A^\times$ not of the form $w=\ao_0^n$,
$$\lim_{z\rightarrow 0} L_w(z) = 0\ .$$
\end{enumerate}
The weight of the function $L_w(z)$ is defined to be the number of
letters which occur in $w$. The functions $L_w(z)$ are defined
inductively by the weight: if $L_w(z)$ has already been defined,
then $L_{\ao_iw }(z)$ is uniquely determined by the differential
equation $(2)$ since the constant of integration is given by $(3)$.
\end{defn}

It follows from the definitions that
$$L_{\ao_i}(z) = \log (z-\sigma_i) -\log (\sigma_i)\ .$$
Later,  we will only  consider linear combinations of functions
$L_w(z)$ which are single-valued on the real interval $(0,\infty)$.
Any such function is uniquely defined on $(0,\infty)$ after having
fixed the branch of the logarithm function $L_{\ao_0}(z)=\log(z)$.

\subsubsection{The shuffle product}
The shuffle product, denoted $\sha$, is a commutative multiplication
law on $\Q\langle A\rangle$ defined inductively by the formulae:
$$w \sha e  = e\sha w =w\quad \hbox{ for all } w \in A^\times\ ,$$
$$   \ao_i w_1\sha \ao_j w_2 = \ao_i( w_1\sha \ao_j w_2)+\ao_j( \ao_i w_1\sha  w_2)\quad \hbox{ for all } w_1,w_2 \in A^\times, \ao_i,\ao_j\in A\ . $$
 We extend the definition of the
functions $L_w(z)$ to $\Q\langle A\rangle $ by linearity:
$$L_{w} (z) = \sum_{i=1}^m q_i L_{w_i}(z)\quad \hbox{ where } w= \sum_{i=1}^m q_i w_i \ , q_i \in \Q\ ,\ w_i \in A^\times\ .$$
The following lemma is well-known.
\begin{lem}The functions $L_w(z)$ satisfy the shuffle relations:
$$L_{w_1}(z) L_{w_2}(z) = L_{w_1 \sha w_2}(z)\ , \hbox{ for all } w_1, w_2 \in \Q\langle A
\rangle \ , \hbox{ and }
  z \in \C\backslash \Sigma\ .$$
\end{lem}
\begin{defn} Now let us define
$$\Or_\Sigma = \Q\Big[ z, {1 \over z-\sigma_0},\ldots, {1 \over z-\sigma_N}\Big]\ ,$$
and let  $L(\Sigma)$  be the $\Or_\Sigma$-module generated by the
functions $L_w(z)$, for $w\in A^\times$.
\end{defn} The shuffle product makes $L(\Sigma)$ into a commutative
algebra. It is a differential algebra for the
 operator $\partial/\partial z$, and is graded by the   weight.

\begin{thm}
The functions $\{L_w(z)\}_{w\in A^\times}$,  are linearly
independent over $\C\otimes_\Q\Or_\Sigma$.
\end{thm}
It follows from the theorem that $L(\Sigma)$ is a polynomial ring,
and a convenient  polynomial basis is given by the functions
$L_w(z)$, where $w$ are  Lyndon words \cite{Re}.

 In order to take primitives in $L(\Sigma)$, we have to enlarge the ring of coefficients slightly. Therefore we must consider
$$\Or_{\Sigma}^+ = \Or_\Sigma \Big[\sigma_{i}, {1 \over \sigma_i -\sigma_j} \
, \ 0\leq i< j\leq N\Big]\ ,$$ and let $L^+(\Sigma)=\Or^+_\Sigma
\otimes_{\Or_\Sigma} L(\Sigma)$ be $\Or_\Sigma^+$-module spanned by
$L_{w}(z)$, for $w\in A^\times$.
\begin{thm} \label{thmhyperlogprim}  Every element $f$ in $L(\Sigma)$ of
weight $n$ has a primitive in $ L^+(\Sigma)$  of weight at most
$n+1$, i.e., an element $F\in L^+(\Sigma)$ such that
$$ \partial F/\partial z= f \ .$$
\end{thm}

The primitive of any generator $f(z) L_w(z)\in L(\Sigma)$, where
$f(z)\in \Or_\Sigma$, can be found explicitly by decomposing $f(z)$
into partial fractions in $\Or_\Sigma^+$. The fact that one must
enlarge $\Or_{\Sigma}$ to $\Or_{\Sigma}^+$ is clear from the
identity:
$${1 \over (z-\sigma_i)(z-\sigma_j)}= {1 \over \sigma_i-\sigma_j} \Big({1 \over z-\sigma_i} - {1 \over z-\sigma_j} \Big)\ . $$
 One is thereby
reduced to the case of finding a primitive of functions of the form
$$(z-\sigma_i)^n L_w(z) \ , \hbox{ where } n \in \Z\ .$$
In the case $n=-1$, a primitive is given by $L_{\ao_iw}(z)$ by
definition $(2)$ above. In all other cases, integration by parts
enables one to reduce to the case of finding a primitive of a
function of lower weight. Thus a primitive of $f(z) L_w(z)$ can be
computed algorithmically in at most $n$ steps, where $n$ is the
weight of $L_w(z)$.

\subsubsection{Logarithmic regularization at infinity}
\label{sectlogreg} Having fixed a branch of the logarithm
$L_{\ao_0}(z)$ we can define the regularization of hyperlogarithms
at infinity.

\begin{prop} Every function  $f (z)\in L(\Sigma)$ can be uniquely
written in the form
$$f(z)= \sum_{i=0}^m f_i(z) \log^i(z) \ ,$$
where $f_i(z)$ is holomorphic at $z=\infty$, for $0\leq i \leq m$.
\end{prop}
We can therefore define the regularized value of $f$ at infinity to
be:
$$\Reg_{z=\infty} f(z) = f_0(\infty)\ .$$
Clearly, $\Reg_{z=\infty} L_{\ao_0}(z)=0$, and for all $1\leq i\leq
N$, $\Reg_{z=\infty} L_{\ao_i}(z) = -\log(\sigma_i).$ Note that the
regularization operator $\Reg_{z=\infty}$ respects multiplication.

There is an analogous notion of regularization at $z=0$. Note that
by the definition of the functions $L_w(z)$ (properties $(1)$ and
$(3)$), we always have
\begin{equation}\label{regat0}
\Reg_{z=0} L_w(z) =0 \ , \quad \hbox{ for all } w \in A^\times\ ,\  w \neq e\ .
\end{equation}

\begin{defn} Suppose that $f(z) \in L(\Sigma)$ is
  holomorphic on the real interval $(0,\infty)$, and that $f(z)\,dz$ has at most logarithmic singularities at $z=0,
\infty$.  We define the regularized integral of $f(z)\,dz$ along
$[0,\infty]$  to be:
$$ \int_0^\infty f(z) \,dz = \Reg_{z=\infty} F(z) - \Reg_{z=0} F(z)\ , $$
where $F(z)\in L(\Sigma)$ is a primitive of $f(z)$. The integral
converges. In practice, we write  $F(z)$ in terms of the basis of
functions $L_w(z)$, and choose the constant term  to be zero. It
will follow that $\Reg_{z=0} F(z)$  vanishes, and the integral is
simply given by $\Reg_{z=\infty}F(z)$ in this case.
\end{defn}

It is clear that the regularized integral is additive. The point is
that, in order to compute an integral which is convergent, one is
allowed to break it into a sum of logarithmically divergent pieces,
compute the regularized integral of each, and add the answers
together. This is illustrated in $\S\ref{sect3spokes}$.

\subsubsection{Multiple Zeta Values}\label{sectMZVS} Consider the  case  $N=1$,
and   $\sigma_0=0$,  $\sigma_1=-1$. After a change of variables
$z\mapsto -z$,  we retrieve the classical situation on
$\Pro^1\backslash\{0,1,\infty\}$. The functions $L_w(-z)$ are
therefore multiple polylogarithms in one variable \cite{Br2}, and we
can compute the regularized values at infinity in terms of multiple
zeta values using standard methods. To do this, let
$X=\{\x_0,\x_1\}$ denote an alphabet with two letters, and let
$Z(\x_0,\x_1)$ be Drinfeld's associator:
$$Z(\x_0,\x_1)= \sum_{w\in X^\times} \zeta_{\sha}(w)\, w \ ,$$
where the numbers $\zeta_{\sha}(w)\in \R$ satisfy:
$$\zeta_{\sha}(\x_0) = \zeta_{\sha}(\x_1)= 0\ ,\quad \hbox{ and }
 \quad \zeta_{\sha}(\x_0^{n_k-1}\x_1\ldots
\x_0^{n_1-1} \x_1) = \zeta(n_1,\ldots, n_k)\ ,$$ for all integers
$n_1,\ldots, n_k\geq 1$ such that $n_k\geq 2$. Extending
$\zeta_{\sha}(w)$ by linearity to $\Q\langle X \rangle$, we have the
shuffle relations
$$ \zeta_{\sha} (w_1)\zeta_{\sha}(w_2) = \zeta_{\sha}(w_1\sha w_2)\quad \hbox{ for all } w_1,w_2 \in X^\times\
.$$ The numbers $\zeta_{\sha}(w)$ are uniquely determined by these
properties. Using the fact that $\zeta(1,2)=\zeta(3)$, one can
verify that, up to weight three:
$$Z(\x_0,\x_1)= 1+ \zeta(2)\, (\x_0\x_1 - \x_1\x_0) + \zeta(3)\,(\x_0^2\x_1 + \x_0\x_1^2 + \x_1\x_0^2+\x_1^2\x_0-2\x_0\x_1\x_0 - 2\x_1\x_0\x_1 ) +\ldots $$

\begin{lem} \label{lemzinf} Let $\Sigma=\{0,-1,\infty\}$, and let $w\in X^\times$.
Set
\begin{equation}\label{zinfdef}
\zeta_\infty(w) = \Reg_{z=\infty}L_w(z) \ .\end{equation} Then the
generating series of regularized values at infinity is given by:
$$ \sum_{w\in X^\times} \zeta_{\infty}(w)\, w =Z^{-1}
(\x_1 -\x_0,-\x_0)\ .$$
\end{lem}
A simple computation shows that:
\begin{equation}\label{AssocZatinf}
Z^{-1}(\x_1-\x_0,-\x_0) = 1 + \zeta(2)\, (\x_1\x_0-\x_0\x_1) +
\zeta(3) \, (\x_0\x_1^2 +\x_1^2\x_0 - 2\, \x_1\x_0\x_1)+ \ldots
\end{equation}
which allows us to read off the values of $\zeta_\infty(w)$ for $w$
up to weight 3.
\newpage

\subsection{An algebra of polylogarithms for reducible graphs}
The integration process takes place in an algebra of polylogarithm
functions in several variables. A basis of such  functions will be
products of hyperlogarithms $L_{w_1}(z_1)\ldots L_{w_k}(z_k)$, each
of which is viewed as a function of a single variable $z_i$, where
$z_i$ is some Feynman parameter. Now, however, the singularities
$\sigma_0,\ldots, \sigma_N$ will themselves be rational functions in
the remaining Feynman parameters $\alpha_1,\ldots,
\alpha_k$.\footnote{We use the word \emph{hyperlogarithm} to denote
a function $L_w(z)$, considered as a function of the single variable
$z$, with  constant singularities $\sigma_i\in \C$. In contrast, we
use the word \emph{polylogarithm} to denote $L_w(z)$ also, but this
time viewed as a function of several variables $z$ and
$\alpha_k,\ldots,\alpha_N$ where some of the singularities
$\sigma_i$ depend on the $\alpha_k,\ldots, \alpha_N$. }

 Let $S=\{f_1,\ldots, f_M\}$  denote  a Fubini reducible  set of
polynomials in the variables $\alpha_1,\ldots, \alpha_N$. We can
assume that the algorithm of $\S\ref{sectStrongRed}$   terminates
when we reduce $S$ with respect to  the variables
$\alpha_{1},\ldots, \alpha_{N}$ in that order. We define a nested
sequence of rings $R_N \subset R_{N-1} \subset \ldots \subset R_1$
recursively, as follows. Set $R_{N+1}=\Q$ and define:
\begin{equation} R_{k} = R_{k+1} \big[\alpha_{k}, {1\over
\alpha_{k}}, {1 \over f }: f\in S_{[\alpha_1,\alpha_2,\ldots,
\alpha_k]} \big]\ , \quad \hbox{for } 2\leq k\leq N\ .
\end{equation}
 The  fibration algorithm
ensures that  every $f \in S_{[\alpha_1,\alpha_2,\ldots, \alpha_k]}$
is linear in $\alpha_k$. Thus, if $\{f_1,\ldots, f_{m_k}\}$ is the
set of  elements in $S_{[\alpha_1,\alpha_2,\ldots, \alpha_k]}$ whose
coefficient of  $\alpha_k$ is non-zero, we can write:
$$f_i= a_i \,\alpha_k + b_i\ ,\hbox{ where } a_i,b_i \hbox{ are invertible in } R_{k+1}\ ,$$
for all $1\leq i\leq m_k$. As in $\S\ref{sectStrongRed}$, we set
\begin{equation}
\Sigma_k= \{\sigma_1,\ldots, \sigma_{m_k} \} \ , \hbox{ where }
\sigma_i = -{b_i \over a_i}\ , \end{equation} and with these
notations, we can write: \begin{equation} R_{k} = R_{k+1}
\big[\alpha_{k}, {1\over \alpha_{k}-\sigma_i},\, \sigma_i \in
\Sigma_k \big]
\end{equation}
We refer the reader to $\S\ref{sect3spokesreduc}$ for a  worked
example in the case of the wheel with 3 spokes diagram.

\begin{defn} Let $L(R_N) = R_N$, and define inductively
\begin{equation}L(R_k) = L(\Sigma_k) \otimes L(R_{k+1})\quad \hbox{ for } 1\leq
k\leq N\ .\end{equation}
\end{defn}

Every element of $L(R_k)$ can be represented as a sum of terms of
the form:
\begin{equation}
\phi= f(\alpha_k,\ldots, \alpha_N) \, L_{w_1}(\alpha_k)
L_{w_2}(\alpha_{k+1}) \ldots L_{w_{N-k+1}} (\alpha_N)\
,\end{equation} where $f\in R_{k}$ is a rational function of
$\alpha_k,\ldots, \alpha_{N}$, and $L_{w_i}(\alpha_{k+i-1})\in
L(\Sigma_{k-i+1})$ is a hyperlogarithm in $\alpha_{k+i-1}$, whose
set of  singularities $\Sigma_{k-i+1}$ are rational functions of the
higher variables $\alpha_{k+i},\ldots, \alpha_N$.

\begin{defn} The weight of the element $\phi$ is
$|w_1|+\ldots+ |w_N|$, where $|w|$ denotes the number of
letters in a word $w$.
\end{defn}

Every element $\phi\in L(R_k)$ is a multivalued function on a
certain open subset of $\C^{N-k}$.
When we integrate, we will only consider elements $\phi$ which are
holomorphic, and hence single-valued, on the real hypercube
$(0,\infty)^{N-k}\subset \C^{N-k}$.

In this manner, we have defined an explicit  algebra of
polylogarithm functions associated to a reducible Feynman graph $G$.
Its elements can be represented symbolically as a linear combination
of products of words. We next show how to compute any  integral of
the type $(\ref{Intrologint})$ associated to $G$ by working inside this
algebra.

\subsection{Existence of primitives}
The integration process requires finding primitives at each stage.
As explained earlier, we need to enlarge the coefficients of the
algebra $L(R_k)$ slightly in order to do this. In the above
notation, we set
$$R_k^+ = R_k \Big[ \sigma_i,\, {1 \over \sigma_i -\sigma_j}\ , \ \hbox{for }\sigma_i,\sigma_j \in \Sigma_k \hbox{ distinct}\Big]\ .$$
The following theorem establishes the existence of primitives.
\begin{thm} \label{thmpolylogprim} Let $f\in L(R_k)$ denote a function of weight $n$ in the
variables $\alpha_k,\ldots, \alpha_{N}$. Then there exists a unique
function $F\in R^+_k \otimes_{R_k}L(R_k)$ of weight at most $n+1$,
which is a primitive of $f$ with respect to $\alpha_k$, and is
regularized at $0$:
$${ \partial F\over \partial \alpha_k}  = f\ , \ \hbox{ and }\quad  \Reg_{\alpha_k=0} F = 0\ .$$
\end{thm}
The construction of the primitive follows  immediately from theorem
$\ref{thmhyperlogprim}$. We can assume $f$ is a generator of
$L(\Sigma_k)$ of the form $f=f_1 \, g$,
 where $f_1\in L(\Sigma_k)$,  and $g\in L(R_{k+1})$ is a function of $\alpha_{k+1},\ldots,\alpha_{N}$ only.
Theorem $\ref{thmhyperlogprim}$ provides a primitive $F_1\in R^+_k
\otimes_{R_k}L(R_k)$ for $f_1$. The required primitive of $f$ is
therefore $F=F_1\,g$.

\subsection{Restricted Regularization}\label{sectRReg}
In the notations of theorem $\ref{thmpolylogprim}$, we must next
 consider the regularized limit of such a primitive $F$ at
$\alpha_k = \infty$,
$$\int_{0}^{\infty} f\, d\alpha_k =  \Reg_{\alpha_k=\infty} F\ .$$
 This is a function of $\alpha_{k+1},\ldots, \alpha_{N}$, but will not  lie in $\C(\alpha_{k+1},\ldots, \alpha_N)\otimes_\Q L(R_{k+1})$, since it will have
extra singularities corresponding to the loci $\sigma_i=\sigma_j$,
where $\sigma_i,\sigma_j \in\Sigma_k$.
 The main technical point is
that  these extra singularities will cancel during the iterated
integration process. Therefore, by mapping superfluous singular
terms to zero, we can compute the integrals by using a more
economical space of functions.

In this way, we will define a restricted regularization map:
\begin{equation}\label{defnRReg}
\RReg_{\alpha_k=\infty} : L(R_k) \To \C(\alpha_{k+1},\ldots,
\alpha_N)\otimes_\Q L(R_{k+1})\ ,\end{equation} which consists of
first taking  the regularization at infinity  as defined in
$\S\ref{sectlogreg}$ and then projecting onto the algebra
$\C(\alpha_{k+1},\ldots, \alpha_N) \otimes_\Q L(R_{k+1})$.

The projection map is easily illustrated for functions of a single
variable, i.e., in the case of hyperlogarithms. Consider two sets of
distinct points in $\C$:
$$\Sigma=\{\sigma_0,\sigma_1,\ldots, \sigma_n\} \subset \Sigma'=\{\sigma_0,\ldots, \sigma_n, \sigma_{n+1},\ldots, \sigma_m\}\ ,$$
and let $L(\Sigma), L(\Sigma')$ denote the corresponding
hyperlogarithm algebras, indexed by the set of words in the
alphabets $A=\{\ao_0,\ldots, \ao_n\}\subset A'=\{\ao_0,\ldots,
\ao_m\}$, respectively. There is a projection map
\begin{eqnarray}
\pi_{\Sigma}:L(\Sigma')& \To& \Or_{\Sigma'}\otimes_\Q L(\Sigma) \\
L_{w}(z) & \mapsto & L_{\pi_A(w)}(z)\ , \nonumber
\end{eqnarray}
where $\pi_A(w)=w$ if $w\in A^{\times}$, but $\pi_A(w)=0$ if $w$
contains a letter  in $A'\backslash A$. The map $\pi_{\Sigma}$ is a
homomorphism for the shuffle product.  Thus, for any set  $\Sigma'$
containing $\Sigma$, we can in this way project out any
singularities in $\Sigma'$ which are not in $\Sigma$.

The idea is to apply this argument to $\Reg_{\alpha_k=\infty} F$. If
one considers $\alpha_{k+2},\ldots,\alpha_{N}$ to be constant, then
one can verify that $\Reg_{\alpha_k=\infty} F$ is a hyperlogarithm function in the variable
$\alpha_{k+1}$. We can therefore write it in the form:
$$\Reg_{\alpha_k=\infty} F=\sum_i q_i\, L_{w_i}(\alpha_{k+1})\ , \hbox{ where }  L_{w_i}(\alpha_{k+1})\in L(\Sigma')\ ,$$
where $L_{w_i}(\alpha_{k+1})$ has singularities in some set
$\Sigma'$ which may be larger than $\Sigma_{k+1}$. After applying
the projection map $\pi_{\Sigma_{k+1}}$, we have
$$\pi_{\Sigma_{k+1}}\big(\Reg_{\alpha_k=\infty} F\big)=\sum_i q_i\, L_{w_i}(\alpha_{k+1})\ , \hbox{ where }  L_{w_i}(\alpha_{k+1}) \in L(\Sigma_{k+1})\ .$$
The coefficients $q_i$ are  functions of $\alpha_{k+2},\ldots,
\alpha_{N}$, and can themselves be projected down to
$L(\Sigma_{k+2})$, \ldots, $L(\Sigma_N)$ in turn. Thus we can
define:
\begin{equation}\label{defRReg}
\RReg_{\alpha_k=\infty} F = \pi_{\Sigma_N} \circ \ldots \circ
\pi_{\Sigma_{k+2}}\circ  \pi_{\Sigma_{k+1}} \big(
\Reg_{\alpha_k=\infty} F \big)\ .
\end{equation}
Therefore  $\RReg_{\alpha_k=\infty} F$  lies in
$\C(\alpha_{k+2},\ldots,\alpha_N)\otimes L(R_{k+1}).$

\begin{thm}
If $S$ is unramified (definition $\ref{defnunram}$), then the
coefficients of $\RReg_{\alpha_k=\infty} F$ are multiple zeta
values.
\end{thm}

The restricted regularization $\RReg_{\alpha_{k=\infty}}\, f$ for
any element $f\in L(R_{k})$ can be computed algorithmically by
succesive differentation with respect to $\alpha_{k+1},\ldots,
\alpha_N$.  This uses an induction on the weight which is
illustrated on some examples in $\S\ref{sectEXreg}$ below. The only
difficulty in practice is to compute the constant terms
$$\Reg_{\alpha_N=0}\circ\ldots\circ
\Reg_{\alpha_{k+1}=0}\Big(\Reg_{\alpha_k=\infty} F\Big)\in \C\ .$$
The statement of the theorem is that this number lies in $\MZV$
under suitable conditions. The proof uses associators in higher
dimensions generalizing the ideas behind $\S\ref{sectMZVS}$.





\subsection{The integration algorithm}
We wish to compute
$$I = \int_{0}^{\infty}\cdots\int_{0}^{\infty}  f_1 \,d\alpha_1\ldots d\alpha_N\ , \hbox{ where } f_1 \in R_1\ .$$
The idea is to integrate, one variable at a time, using the
polylogarithm functions in $L(R_1)$. At the $k^{\mathrm{th}}$ stage
of the integration process, we will have:
$$I = \int_{0}^{\infty}\cdots\int_{0}^{\infty} f_k \, d\alpha_k\ldots d\alpha_N\ , \hbox{
where } f_k \in \MZV \otimes_\Q L(R_k)\ ,$$ and $f_k$ has weight at
most $k$. The integrands $f_k$ are calculated recursively as
follows:
\begin{enumerate}
  \item The function $f_k$ has a primitive $F_k \in R^+_k\otimes_{R_k}
L(R_k)$ of weight at most $k+1$.
  \item Since we can assume that $\Reg_{\alpha_k=0} F_k=0$, we have
$$\int_{0}^\infty f_k \,d\alpha_k = \Reg_{\alpha_k=\infty} F_k\ .$$
 In general, however, $\Reg_{\alpha_k=\infty} F_k$ is not an element
of $\MZV\otimes_\Q L(R_{k+1})$, so we must use restricted
regularization instead.

 \item Therefore, we   define
$$f_{k+1} = \RReg_{\alpha_k=\infty} F_k \in  \MZV  \otimes_\Q L(R_{k+1})\ .$$
Although $\int_{0}^{\infty} f_k \,d{\alpha_k} $ is not in general
equal to $f_{k+1}$, one can prove that
$$ \int_{0}^{\infty}\cdots\int_{0}^{\infty} f_k \, d\alpha_k\ldots d\alpha_N =  \int_{0}^{\infty}\cdots\int_{0}^{\infty} f_{k+1} \, d\alpha_{k+1}\ldots d\alpha_N \ ,$$
and the induction goes through.
\end{enumerate}

At the last stage of the integration, we  have
$$I=\int_{0}^{\infty} f_N \,d\alpha_N= \Reg_{\alpha_N=\infty} F_N\ ,$$
where $F_N$ is the primitive of $f_N$ defined in $(2)$ above. Since
$F_N\in \MZV\otimes L(R_N)$, its regularized value at $\infty$ is a
$\Q$-linear combination of multiple zeta values by
$\S\ref{sectMZVS}$. This proves that $I\in \MZV$.

Since the weight of the integrand increases by at most one at each
integration, we can use the above argument to bound the weights of the periods we obtain.

\begin{rem} There is never any need to choose branches of any hyperlogarithm functions, because one can
prove that every function $f_k$ which occurs in the above process is
in fact  holomorphic on the open hypercube
$$\{(\alpha_k,\ldots, \alpha_N)\in \R^{N-k} : \, 0< \alpha_k,\ldots, \alpha_N<\infty\}\ .$$
The functions $f_k$ can, however, have logarithmic singularities
along the codimension 1 faces of this hypercube. A detailed study of
the locus of singularities in the Feynman case, and the
corresponding compactifications, will be given in \cite{Br4}.
\end{rem}

Note that the entire integration process is algorithmic, from the
calculation of primitives to the regularizations at infinity, and
can be reduced to a sequence of elementary manipulations on symbols.

\newpage
\section{Examples of hyperlogarithms  and regularization}
\label{sectEXAMPLES}

We  consider some one-dimensional examples which will occur in the
calculation of the wheel with three spokes diagram. First let $N=1$,
$\sigma_0=0$, and $\sigma_{1}=-1$, and denote the corresponding
alphabet by $X=\{\x_0, \x_1\}$. We set:
\begin{equation}\label{exdefsigmax}
\Sigma_{x}=\{0,-1\} \ , \hbox{ and }    \ \Or_{\Sigma_x} = \Q\Big[
x, {1 \over x}, {1 \over x+1}\Big]\ ,
\end{equation}
 and $L(\Sigma_x)$ is the
$\Or_{\Sigma_x}$-module spanned by the hyperlogarithms $\Li_w(x)$,
where $w$ is any word in $X$. Thus, in $L(\Sigma_x)$, there are
exactly two hyperlogarithms of weight $1$:
$$\Li_{\x_0}=\log(x)\quad  , \quad \Li_{\x_1}=\log(x+1) \ .$$
Note the absence of a minus sign in front of $\Li_{\x_1}$, which
differs  from the usual convention $\Li_1(z)=-\log(1-z)$. In weight two there are
precisely four functions:
$$\Li_{\x_0^2}(x)  \ ,  \Li_{\x_0\x_1}(x)  \ , \  \Li_{\x_1\x_0}(x)\ , \ \Li_{\x_1^2}(x) \ .$$
Using  definition $\ref{defnHyperlog}$, one verifies that these can
also be expressed as:
 $$ {1\over 2} \log^2(x) \ , \ -\Li_2(-x) \ , \   -\Li_2(x+1) \ , \ {1 \over 2 } \log^2(x+1)\ ,$$
respectively. 
They are related by a   single shuffle product:
$$\Li_{\x_0}(x) \Li_{\x_1}(x) = \Li_{\x_0\x_1}(x) +
\Li_{\x_1\x_0}(x)\ .$$ Finally, in weight three there are exactly 8
linearly independent hyperlogarithms. One can check, for example,
that $\Li_{\x_0^2\x_1}(x)=-\Li_3(-x)$.
 If we denote the regularized values at infinity by $\zeta_{\infty} (w) =
\Reg_{x=\infty} L_w(x)$, we have from $\S\ref{sectMZVS}$,
$$\zeta_{\infty}(\x_0^2) = 0 \ , \ \zeta_{\infty}(\x_0\x_1) =
-\zeta(2)  \ , \ \zeta_{\infty}(\x_1\x_0) = \zeta(2)  \ , \
\zeta_{\infty}(\x_0\x_1) = 0\ .$$ Likewise, in weight three, we
deduce from lemma $\ref{lemzinf}$ that:
\begin{equation} \label{Exregzetavalues}
\zeta_{\infty}(\x_1\x_0\x_1)= - 2\, \zeta(3)\ , \
\zeta_{\infty}(\x_0\x^2_1)= \zeta(3)\ , \
\zeta_{\infty}(\x_1^2\x_0)=\zeta(3)
\end{equation}
$$\zeta_{\infty}(\x_0^2\x_1)= 0 \ ,
\zeta_\infty(\x_0\x_1\x_0)=0\ , \  \zeta_\infty(\x_1\x^2_0)=0 \ , \
\zeta_{\infty}(\x_0^3)=0 \ , \ \zeta_{\infty}(\x_1^3)=0\ .$$
\subsubsection{The enlarged case $\{0,\pm 1\}$} Now consider the case when $N=2$, $\sigma_0=0$,
$\sigma_1=-1$, $\sigma_2=1$. We  set
\begin{equation}\label{defnsigmaxplus}
\Sigma_x^+=\{0,-1,1\}    \ , \hbox{ and }    \  \Or_{\Sigma_x^+} =
\Q\Big[x, {1 \over x}, {1 \over x-1}, {1 \over x+1}\Big]\ .
\end{equation}
 We denote the corresponding alphabet by
$X^+=\{\x_0,\x_1,\x_2\}$. The algebra $L(\Sigma^+_x)$ now
  contains 3 hyperlogarithms in weight
1, namely:
 $$\Li_{\x_0}(x)= \log (x) \ , \ \Li_{\x_1}(x)=\log(x+1) \ , \  \Li_{\x_2}(x) =
\log (x-1)\ ,$$ and $3^n$ hyperlogarithms in weight $n$. In this
case, the regularized values at infinity are no longer multiple zeta
values $(\MZV^1)$, but the larger set of alternating multiple zeta values $(\MZV^2)$. Via
the inclusion $L(\Sigma_x) \rightarrow L(\Sigma_x^+)$, we can
identify $L(\Sigma_x)$ as a sub-algebra of $L(\Sigma_x^+)$. In
$\S\ref{sectRReg}$ we defined a projection map:
\begin{eqnarray}
\pi_{\Sigma_x}: L(\Sigma_x^+) &\To&
\Or_{\Sigma_x^+}\otimes_{\Or_{\Sigma_x}}L(\Sigma_x)
\nonumber \\
f(x)\, L_{w}(x) & \mapsto & f(x) \, L_{w}(x) \quad \hbox{ if } w \in
X^\times \ , \nonumber \\
f(x)\, L_{w}(x) & \mapsto & 0 \quad \hbox{ if } w \hbox{ contains
the letter } \x_2 \ . \nonumber
\end{eqnarray}

\newpage
\subsection{Second example} \label{sectexofregs} We now consider an example in two variables. Let
\begin{equation}\label{defExsigmaY}
\Sigma_y=\{\sigma_0,\sigma_1,\sigma_2,\sigma_3\} \  , \hbox{ where }
 \sigma_0=0\ , \ \sigma_1=-1\ ,\ \sigma_2=-x\ ,\
\sigma_3=-{x\over x+1}\ .
\end{equation}
 We set
$$\Or_{\Sigma_y} = \Q\Big[ y, {1 \over y}, {1 \over y+1}, {1 \over
x+y}, {1 \over y+{x\over x+1}}\Big]\ .$$ Let us denote the
corresponding alphabet by $Y=\{\y_0,\y_1,\y_2,\y_3\}$, and thus
$$L_{\y_0}(y) = \log(y)\ , \ L_{\y_1}(y)= \log(y+1) \ , \ L_{y_2}(y) =\log(x+y)-\log(x) \ , $$
$$ \ L_{y_3}(y) = \log( xy+x+y) -\log(x)\ . $$
Thus $L(\Sigma_y)$ is spanned (as a vector space) in weight $n$ by
exactly $4^n$ hyperlogarithm functions $L_w(y)$, which can be
considered as functions of the single variable $y$, for constant
$x$. If we wish to consider the dependence on $x$, we must remove
the singular locus where the $\sigma_i$ collide.

Recall that $\Or_{\Sigma_y}^+$ was defined to be
$\Or_{\Sigma_y}[(\sigma_i), (\sigma_i-\sigma_j)^{-1}]$, and so
\begin{equation}\label{defexosplus}
\Or_{\Sigma_y}^+ = \Q\Big[ x, {1 \over x}, {1 \over x+1} , {1 \over
x-1}\Big] \Big[ y, {1 \over y}, {1 \over y+1}, {1 \over x+y}, {1
\over y+{x\over x+1}}\Big]\ . \end{equation} This is geometrically a
linear fibration over $\Or_{\Sigma^+_x}$ and is pictured below:
\begin{figure}[h!]
  \begin{center}
    \epsfxsize=7.0cm \epsfbox{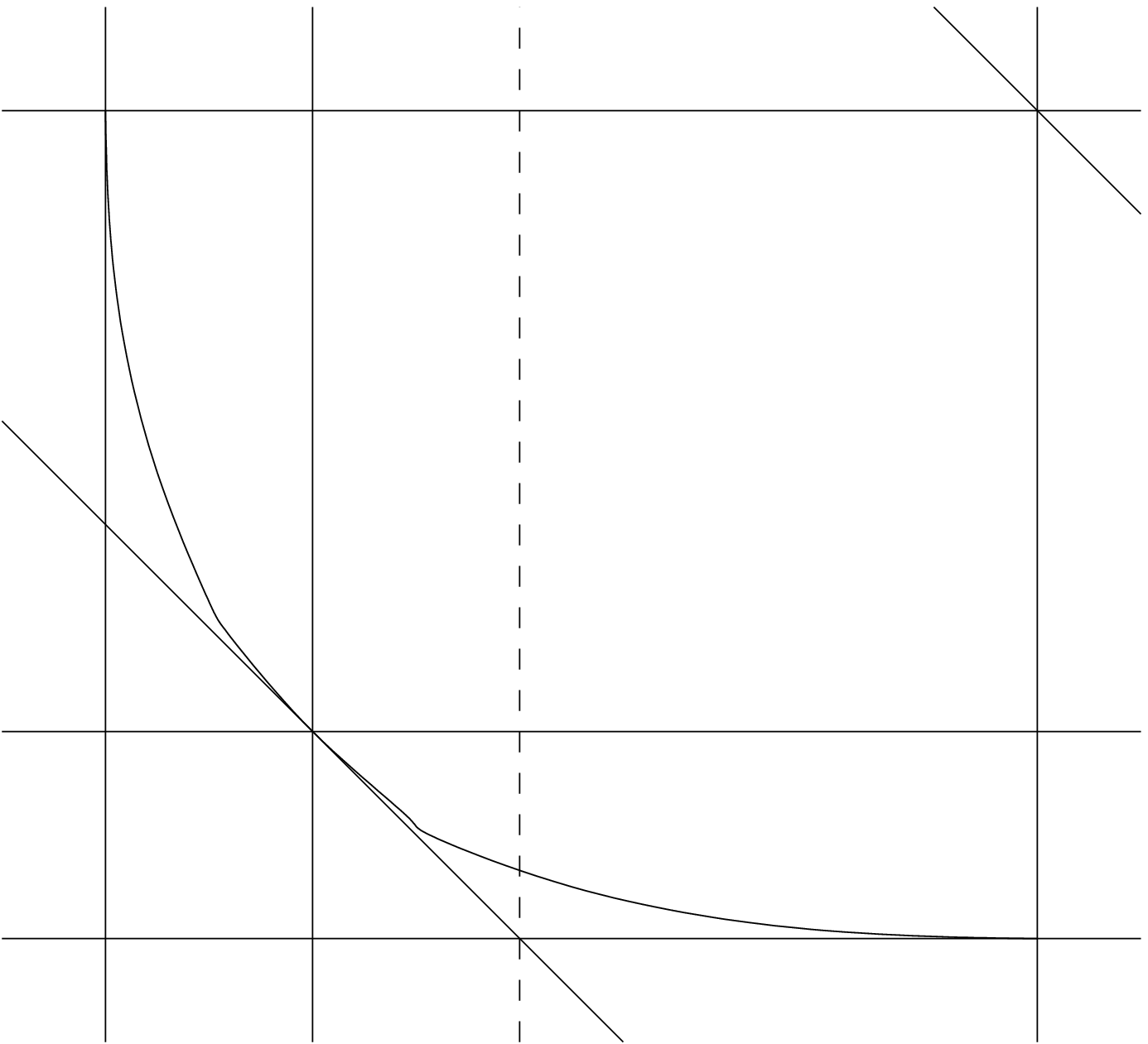}
  \label{figurefib}
\put(-240,50){$y=0$} \put(-245,20){$y=-1$}
\put(-240,160){$y=\infty$}
 \put(-160,-8){$x=0$}
\put(-197,-8){$x=-1$} \put(-35,-8){$x=\infty$}
  \end{center}
\end{figure}

Let $U=\C^2\backslash\{x=0,x=\pm 1,y=0,y=-1, x+y=0, xy+x+y=0 \}$, as
shown.  The elements of $L(\Sigma_y)$ can be viewed as multi-valued
functions on $U$. The full space of polylogarithms is a product of
two hyperlogarithm algebras:
$$L(\Or_{\Sigma_y}^+) = L(\Sigma_{x}^+) L(\Sigma_{y})\ .$$
It is generated by elements $f(x,y) \, L_{w_1}(x) L_{w_2}(y)$ where
$w_1 \in X^\times$, $w_2\in Y^\times$, and $f(x,y) \in
\Or_{\Sigma_y}^+$. The set of elements of weight two are:
$$\Li_{\x_i\x_j}(x) \ , \ \Li_{\x_i}(x) \Li_{\y_i}(y) \ , \
\Li_{\y_i\y_j}( y) \ ,$$ where the $\x$'s are in $X^+$, and the
$\y$'s are in $Y$. There are, respectively, $9, 12,$ and  $16$ such
elements, giving a total of $37$ polylogarithms of weight 2.


\subsubsection{Regularized values}\label{sectEXreg}
Let $\MZV$ denote the $\Q$-algebra spanned by all multiple zeta
values. We defined a regularization map
$$\Reg_{y=\infty} :L(\Sigma_y) \To \MZV \otimes_\Q L(\Or_{\Sigma_y}^+)\ .$$
The regularized values at $y=0$ of all functions $L_w(y)$, $w\in
Y^\times$ vanish. At $y=\infty$,
\begin{equation}\label{exwt1reg} \Reg_{y=\infty} L_{y_0}(y)=0\ , \
\Reg_{y=\infty} L_{y_1}(y)=0\ , \end{equation}
$$\Reg_{y=\infty} L_{y_2}(y)=-\log(x) \ , \ \Reg_{y=\infty} L_{y_3}(y)=\log(x+1)-\log(x)\ .$$
We now compute the regularised values at infinity of some functions
of weight 2.
 To compute $\Reg_{y=\infty} L_w(y)$, we can differentiate with
respect to $x$, which gives a function of lower weight, compute the
regularised value at infinity of this function by induction, and
take the primitive with respect to $x$. In other words,
$$\Reg_{y=\infty} L_w(y) = \int  \Reg_{y=\infty} \Big({\partial \over \partial x} L_w(y)\Big)\, dx\ .$$
The constant of integration is determined from the regularized
values:
$$ \Reg_{y=\infty} \Reg_{x=0} L_w(y) \in \MZV \ ,$$
which can be  calculated by  a generalization of the associator
argument
of $\S\ref{sectMZVS}$. 
\subsubsection{First example}
Let $w =\y_0\y_2$. To compute ${\partial\over \partial x} L_{w}(y)$, we have:
$${\partial \over \partial y} {\partial \over \partial x} L_{\y_0\y_2}(y) = {\partial \over \partial x} {\partial \over \partial y} L_{\y_0\y_2}(y) =
{\partial \over \partial x} {1 \over y} L_{\y_2}(y) = {1 \over y} \Big({1\over x+y} - {1 \over x}\Big)\ . $$
Thus
$${\partial \over \partial x} L_{\y_0\y_2}(y) = \int {-1 \over x(x+y)} \, dy = -{1 \over x} L_{\y_2}(y)\ .$$
The constant of integration is determined by the fact that $L_{\y_0\y_2}(y)$, and hence ${\partial \over \partial x} L_{\y_0\y_2}(y)$, vanishes along $y=0$.
Taking the regularized value at $y=\infty$ gives
$$\Reg_{y=\infty} \Big({\partial \over \partial x} L_{\y_0\y_2}(y)\Big) = -{1 \over x} \Reg_{y=\infty} L_{\y_2}(y) = {\log(x) \over x} \ .$$
By taking a primitive with respect to $x$ (working now in
$L(\Sigma_x^+)$), we deduce that
$$\Reg_{y=\infty} L_{\y_0\y_2}(y) = \int {\log x \over x} dx = L_{\x_0^2}(x)+ c= {1 \over 2} \log(x)^2 +c\ ,$$
where $c$ satisfies $c= \Reg_{x=0} \Reg_{y=\infty} L_{\y_0\y_2}(y)=
0.$
%
It follows that \begin{equation}\Reg_{y=\infty} L_{\y_0\y_2}(y) =
L_{\x_0^2}(x)\ . \end{equation}
\subsubsection{Second example}
Let $w =\y_3\y_1$. By a similar calculation, we have
$${\partial \over \partial x} L_{\y_3\y_1}(y) = {y+1 \over xy+x+y} L_{\y_1}(y)  - {1 \over x+1}L_{\y_3}(y)\ .$$
Taking the regularized limit at $y=\infty$ gives
$${\partial \over \partial x} \Reg_{y=\infty}L_{\y_3\y_1}(y) = {1 \over x+1} \big(\log(x) - \log(x+1)\big)\ .$$
To determine the constant of integration, we use the fact that
$\Reg_{x=0} \Reg_{y=\infty} L_{\y_3\y_1}(y)=\zeta(2)$.
Taking a primitive with respect to $x$ in the ring $L(\Sigma_x^+)$,
we deduce that
$$\Reg_{y=\infty} L_{\y_3\y_1}(y) = L_{\x_1\x_0}(x) - L_{\x_1^2}(x) + \zeta(2)\ . $$
\subsubsection{Further examples}
In a similar way, one can verify  the following:
\begin{eqnarray}
\Reg_{y=\infty} L_{\y_0\y_3}(y) &= & L_{\x_0^2}(x)+L_{\x_1^2}(x)-L_{\x_0\x_1}(x)-L_{\x_1\x_0}(x)\ ,  \\
\Reg_{y=\infty} L_{\y_3\y_2}(y)&= & L_{\x_0\x_1}(x)-L_{\x_1^2}(x)+L_{\x_0^2}(x)\ , \nonumber
\end{eqnarray}
which will be used in the calculation of the wheel with three spokes diagram overleaf.

\subsubsection{Restricted regularization} All the previous examples happen to lie in
$L(\Sigma_x)$.  In the general case, one obtains an answer with an
extra singularity at $x=1$, i.e., a hyperlogarithm in
$L(\Sigma_{x}^+)$. For example, one can check that:
\begin{eqnarray}
 {\partial
\over \partial x} L_{\y_1\y_2}(y)& =&   { L_{\y_1}(y) \over
x(x-1)}-{L_{\y_2}(y)  \over x-1}\ ,  \\
\Reg_{y=\infty} L_{\y_1\y_2}(y) &= &L_{\x_2\x_0}(x)\ . \nonumber
\end{eqnarray}
The reason this happens is because $L_{\y_1\y_2}(y)$ has a singularity at  $\sigma_1=\sigma_2$.
In this case, the restricted regularization map of
$\S\ref{sectRReg}$ gives
$$\RReg_{y=\infty} L_{\y_1\y_2}(y) = \pi_{\Sigma_x} L_{\x_2\x_0}(x)= 0 \ .$$
Since the regularized value of every hyperlogarithm occurring in the
calculation of the wheel with three spokes diagram overleaf already
lies in $L(\Sigma_x)$, we do not in fact need to use the regularized
restriction in the calculation.

\begin{rem}
In order to compute the values  $\Reg_{x=0} \Reg_{y=\infty}
L_{w}(y)$, one can set
$$t_1=-\sigma_2\ , \ t_2=-\sigma_3\ ,  \ \hbox{
and }\  t_3 = -y\ .$$ The hyperlogarithm $L_w(y)$, viewed as a
function of the independent variables $t_1,t_2, t_3$, is then a
unipotent function  on the moduli space
$$\Mod_{0,6} =\mathrm{Spec\,} \Z \Big[ t_1,t_2,t_3, { 1 \over t_1},{1\over t_2},{1\over t_3}, {1 \over 1-t_1},{1 \over 1-t_2},{1 \over 1-t_3},
 {1 \over t_1-t_2}, {1 \over t_1-t_3}, {1 \over t_2-t_3} \Big]\ . $$
Setting $t_1=x, t_2= {x\over x+1}, t_3=-y$ defines a surface inside
$\Mod_{0,6}$, and taking the limit first as  $y\rightarrow 0$ and
then  $x\rightarrow 0$, corresponds to taking the regularized limit
of $L_w(y)$ at the tangential base-point at $t_1=0,t_2=0,t_3=0$
defined by the sector
$$0\leq t_3\leq t_2 \leq t_1 \leq 1\ .$$
Equivalently, this corresponds to a unique  point on the
compactification $\overline{\Mod}_{0,6}$ which lies in the deepest stratum, and defined over the integers. Likewise, taking the
regularized limit at $x=0$, $y=\infty$ corresponds to taking a limit
at a different point in $\overline{\Mod}_{0,6}$. By the results of
\cite{Br3}, this limit can be described by a generalized associator
and is expressible in terms of multiple zeta values. The general
case is similar and follows from the ramification condition
$(\ref{ramcond})$ (this will be discussed in \cite{Br4}).

\end{rem}

\newpage
\section{The wheel with three spokes}

\subsection{The Fubini reduction algorithm} \label{sect3spokesreduc} Let $G$ be the wheel with
three spokes depicted in fig. 1. We take $\lambda=5$, and set
$\alpha_5=1$ from the outset. We therefore have $S'=\{U_{\widetilde{G}}|_{\alpha_5=1}\}$,
where $U_{\widetilde{G}}=U_G\alpha_6+V_G$, and $U_G$, $V_G$ are given by $(\ref{U2p2L})$. We wish to compute integrals of the form:
\begin{equation}
\label{3spokesintegral} I_G = \int_{\alpha_5=1}
{\prod_{i=1}^6 \log^{m_i}(\alpha_i)\log^{m}(U_{\widetilde{G}}) \over U^2_{\widetilde{G}}
} d\alpha_6 d\alpha_1d\alpha_2 d\alpha_3 d\alpha_4\ ,\quad m_i, m\in \N\ , \end{equation} by
integrating with respect to the variables
$\alpha_6,\alpha_1,\alpha_2,\alpha_3,\alpha_4$, in that order. The Fubini
reduction algorithm gives:
\begin{eqnarray} \nonumber
S'_{[6]}&=&\{ \alpha_1\alpha_2+\alpha_1\alpha_4 +
\alpha_1\alpha_3+\alpha_3\alpha_2+\alpha_3\alpha_4+\alpha_2+\alpha_3+\alpha_4\ ,  \\
&\quad
&\alpha_1\alpha_3\alpha_4+\alpha_1\alpha_2\alpha_4+\alpha_3\alpha_2\alpha_4+
\alpha_1\alpha_3+\alpha_2\alpha_3+\alpha_2\alpha_4+\alpha_1\alpha_4+\alpha_1\alpha_2
 \}\ .\nonumber
\end{eqnarray}
After reducing with respect to $\alpha_1$, we have:
\begin{eqnarray} \nonumber
 S'_{[6,1]}=\{& \!\!\!\!\alpha_2 +\alpha_3+ \alpha_4+\alpha_2\alpha_3+\alpha_3\alpha_4\ ,
\,\,\alpha_2+\alpha_3+\alpha_4\ ,\,\,
  \alpha_3\alpha_4+\alpha_3+\alpha_4 \ ,\,\, &\!\!\!\!   \\
 & \!\!\!\!\nonumber\alpha_2\alpha_4+\alpha_3\alpha_4+\alpha_2+\alpha_3+\alpha_4\ , \,\, \nonumber\alpha_3\alpha_4+\alpha_2+\alpha_3+\alpha_4&\!\!\!\!\}\
\end{eqnarray}
Reducing with respect to the variable $\alpha_2$ gives:
$$
 S'_{[6,1,2]}=\{ \alpha_4+1\ ,\,\,
  \alpha_3+1\ , \,\,
 \alpha_3\alpha_4 + \alpha_3+\alpha_4\
,\,\, \alpha_3+\alpha_4 \,\,
 \}\ .
$$
Finally, by reducing with respect to the variable $\alpha_3$, we
obtain:
$$
 S'_{[6,1,2,3]}=\{ \alpha_4+1\}\ .
$$
Note that at each stage, every polynomial is linear with respect to
every variable. The corresponding sets of singularities
$\Sigma_1,\ldots, \Sigma_4$ are therefore:
\begin{equation} \label{Exsigmas}
  \Sigma_1 = \{0, -{ \alpha_2+\alpha_3+\alpha_4+\alpha_2\alpha_3+\alpha_3\alpha_4\over \alpha_2+\alpha_3+\alpha_4}
, -{ \alpha_2\alpha_3+\alpha_2\alpha_4+\alpha_2\alpha_3\alpha_4\over
\alpha_2+\alpha_3+\alpha_4+\alpha_2\alpha_4+\alpha_3\alpha_4} \}\
,\end{equation}
$$\Sigma_2 = \{0, -(\alpha_3+\alpha_4), -(\alpha_3+\alpha_4+\alpha_3\alpha_4), -{\alpha_3+\alpha_4+\alpha_3\alpha_4 \over \alpha_3+1} ,-{\alpha_3+\alpha_4+\alpha_3\alpha_4 \over \alpha_4+1}\}
\ , $$
$$\Sigma_3 =
\{0,-1,-\alpha_4,-{\alpha_4\over \alpha_4+1}\} \quad , \quad
\Sigma_4 = \{0,-1\} \ .$$
 Taking the limits  as $\alpha_2\rightarrow 0$,
$\alpha_3\rightarrow 0$, $\alpha_4\rightarrow 0$, in that order, we
obtain singularities in $0,-1$ only. Therefore the conditions of
theorem $\ref{thmMAINCRITERION}$ are satisfied.

\begin{cor} Every integral of the form $(\ref{3spokesintegral})$  lies in $\MZV$. 
\end{cor}
We have a nested sequence of rings $\Q \subset R_4 \subset R_3
\subset R_2 \subset R_1,$  defined as follows:
\begin{eqnarray}
R_4&= &\Q[\alpha_4, \alpha_4^{-1}, (\alpha_4+1)^{-1}] \ ,\nonumber \\
R_3&= &R_4[\alpha_3, \alpha_3^{-1},
(\alpha_3+1)^{-1},(\alpha_3+\alpha_4)^{-1},(\alpha_3\alpha_4+\alpha_3+\alpha_4)^{-1}]
\ ,\nonumber \\
 R_2&= &R_3[\alpha_2, \alpha_2^{-1},
(\alpha_2 +\alpha_3+
\alpha_4+\alpha_2\alpha_3+\alpha_3\alpha_4)^{-1} ,
(\alpha_2+\alpha_3+\alpha_4)^{-1}, \nonumber  \\
& & \qquad \qquad
(\alpha_2\alpha_4+\alpha_3\alpha_4+\alpha_2+\alpha_3+\alpha_4)^{-1},
(\alpha_3\alpha_4+\alpha_2+\alpha_3+\alpha_4)^{-1}] \ ,\nonumber \\
R_1&= &R_2[\alpha_1, \alpha_1^{-1},
(\alpha_1\alpha_2+\alpha_1\alpha_4 +
\alpha_3\alpha_1+\alpha_3\alpha_2+\alpha_3\alpha_4+\alpha_2+\alpha_3+\alpha_4)^{-1}]
\ , \nonumber \\
& & \qquad \qquad
(\alpha_3\alpha_1\alpha_4+\alpha_3\alpha_2\alpha_4+\alpha_1\alpha_2\alpha_4+
\alpha_1\alpha_3+\alpha_2\alpha_3+\alpha_2\alpha_4+\alpha_1\alpha_4+\alpha_1\alpha_2)^{-1}]
\ .\nonumber
\end{eqnarray}
One can compute any particular integral $(\ref{3spokesintegral})$ by
working in the algebra of polylogarithms $L(R_1) = L(\Sigma_1)\,
L(\Sigma_2)\, L(\Sigma_3)\, L(\Sigma_4),$ where $L(\Sigma_1),\ldots,
L(\Sigma_4)$
 are hyperlogarithm algebras on $3,5,4,2$ letters
respectively. 

\subsection{Calculation of the leading term for the wheel with 3 spokes}
\label{sect3spokes} We illustrate the method by calculating  in
complete detail the Feynman amplitude of the wheel with three spokes
$\widetilde{G}$, and reprove the result, due originally to
Broadhurst and Kreimer, that it evaluates to $6\,\zeta(3)$.

Let $U=U_G$, $V=V_G$ be as in $(\ref{U2p2L})$. We set  $\lambda=\{5\}$, and   wish to compute $(\ref{IG0})$:
$$I_G= \int_{\{\alpha_5=1\}}{1\over UV} \,{d\alpha_1d\alpha_2d\alpha_3d\alpha_4}\ .$$
From equation $(\ref{Ddefn})$, we have
$D=\alpha_2\alpha_5+
\alpha_4\alpha_5+\alpha_3\alpha_4+\alpha_3\alpha_5$, and
$$U^{(1,2)}= \alpha_3\alpha_4+\alpha_3\alpha_5 + \alpha_4\alpha_5\ ,\  U^{(1)}_2 = \alpha_3+\alpha_5\ ,\ U^{(2)}_1= \alpha_3+\alpha_4 \ ,\ U_{12}= 1\ .$$
$$V^{(1,2)}= 0\ ,\  V^{(1)}_2 = V^{(2)}_1= \alpha_3\alpha_4+\alpha_3\alpha_5+\alpha_4\alpha_5  \ ,\ V_{12}= \alpha_4+\alpha_5\ .$$
Let us set $\alpha_5=1$ once and for all. Using the notation $(\ref{defCurlynotation})$ we have:
\begin{eqnarray}
 \{ U^{(1)}_{2}, U^{(1,2)} | D_2, D^{(2)} \}&=& {(\alpha_3+1)
\log(\alpha_3+1) \over \alpha_3(\alpha_3\alpha_4+\alpha_3+\alpha_4)
} +{ \log(\alpha_3\alpha_4+\alpha_3+\alpha_4)\over
\alpha_3\alpha_4+\alpha_3+\alpha_4 } \ ,\nonumber \\
 \{ U_{12}, U^{(2)}_1 | D_2, D^{(2)} \}&=&
{\log(\alpha_3\alpha_4+\alpha_3+\alpha_4)-\log(\alpha_3+\alpha_4)\over
\alpha_3\alpha_4}+{\log(\alpha_3+\alpha_4)\over
(\alpha_3\alpha_4+\alpha_3+\alpha_4)}\ , \nonumber \\
\{ V^{(1)}_{2}, V^{(1,2)} | D_2, D^{(2)} \} &=&{2 \log
(\alpha_3\alpha_4+\alpha_3+\alpha_4) \over
(\alpha_3\alpha_4+\alpha_3+\alpha_4) }\ , \nonumber \\
\{ V_{12}, V^{(2)}_1 | D_2, D^{(2)} \}&=&
{(\alpha_4+1)\log(\alpha_4+1)\over
\alpha_4(\alpha_3\alpha_4+\alpha_3+\alpha_4)}
+{\log(\alpha_3\alpha_4+\alpha_3+\alpha_4)\over
\alpha_3\alpha_4+\alpha_3+\alpha_4} \ .\nonumber
\end{eqnarray}
By corollary $\ref{corstep3}$, we can skip the first two integration steps and go straight to
\begin{equation}\label{appendixStep3}
I_G= \int_0^\infty \int_0^\infty
{1 \over \alpha_3\alpha_4}
\Big(\log(\alpha_3+1)+\log(\alpha_3+\alpha_4)-
\log(\alpha_3\alpha_4+\alpha_3+\alpha_4)\Big)
\end{equation}
$$+{1\over \alpha_4 (\alpha_3\alpha_4+\alpha_3+\alpha_4)}
\Big((\alpha_4+1)\log(\alpha_4+1)-\log(\alpha_3+1) - \alpha_4
\log(\alpha_3+\alpha_4) \Big)
 d\alpha_3 d\alpha_4\ .$$

One can verify  that the integrand has no polar singularities along
the faces of the domain of integration $X=[0,\infty]\times
[0,\infty]$. It has a pole of order exactly  one along the
hypersurface $\alpha_3\alpha_4+\alpha_3+\alpha_4=0$, which meets the
boundary of $X$ in codimension two.
Even though the integral converges, the simplest way to compute it
is to calculate the regularised integral of each individual term,
which may be at most logarithmically divergent, and add the
contributions together.

We set $\alpha_3=y,\alpha_4=x$. In $\S\ref{sectEXAMPLES}$,  $R_3$
and $R_4$ were called $\Or_{\Sigma_y}$ and $\Or_{\Sigma_x}$,
respectively. In the hyperlogarithm notation of
$\S\ref{sectHyperlogarithms}$, the integrand of $I_G$ is an element
of $L(R_4)= L(\Or_{\Sigma_y})$ of weight 1. We   rewrite $I_G$ as
follows:
$$I_G = \int_{0}^\infty \int_{0}^\infty {1 \over xy} \Big(L_{\y_1}(y)+ L_{\y_2}(y)-L_{\y_3}(y)\Big)+{1\over x} L_{\x_1}(x)\Big({1 \over y+{x\over x+1}}\Big)$$
$$-{1 \over x(x+1)}{1 \over y+{x\over x+1}}\Big( L_{\y_1}(y)  +xL_{\y_2}(y)+xL_{\x_0}(x)\Big)dx dy\ .$$
Using theorem $\ref{thmpolylogprim}$, we can take a primitive with
respect to $y$, which  gives:
$$I_G = \int_{0}^\infty {1 \over x}\Reg_{y=\infty} \Big[ L_{\y_0\y_1}(y)+L_{\y_0\y_2}(y)-L_{\y_0\y_3}(y)+L_{\x_1}(x)L_{\y_3}(y)\Big]$$
$$- {1 \over x(x+1)}\Reg_{y=\infty} \Big[ L_{\y_3\y_1}(y) +xL_{\y_3\y_2}(y)+ xL_{\x_0}(x)L_{\y_3}(y)\Big]dx\ .$$
Using  the shuffle relations, and equation $(\ref{exwt1reg})$,
observe that
$$\Reg_{y=\infty}\, L_{\x_1}(x)L_{\y_3}(y)= L_{\x_1}(x)\big( L_{\x_1}(x) - L_{\x_0}(x)\big)=2\,L_{\x_1^2}(x) -L_{\x_1\x_0}(x)-L_{\x_0\x_1}(x)\ ,$$
and similarly,
$$\Reg_{y=\infty}\,
L_{\x_0}(x)L_{\y_3}(y)=L_{\x_0\x_1}(x)+L_{\x_1\x_0}(x)-2\,L_{\x_0^2}(x)\
.$$ All the remaining regularized limits were given in
$\S\ref{sectexofregs}$.
 Substituting them in gives:
$$I_G=\int_{0}^\infty {1 \over x} \Big[ \zeta(2) + L_{\x_0\x_1}(x) +L_{\x_1\x_0}(x)-L_{\x_1^2}(x)+L_{\x_1}(x)\big( L_{\x_1}(x) - L_{\x_0}(x)\big)\Big]$$
$$ -{1 \over x(x+1)}\Big[  \zeta(2)+ L_{\x_1\x_0}(x) - L_{\x_1^2}(x) + x\big(L_{\x_0\x_1}(x)-L_{\x_1^2}(x)+L_{\x_0^2}(x) +L_{\x_0}(x)\big(L_{\x_1}(x)-L_{\x_0}(x)\big)\big) \Big]\ .$$
After decomposing into partial fractions, and expanding out the
shuffle products, we obtain the fourth step:
\begin{equation} \label{fourthstep}
I_G=\int_{0}^\infty {1 \over x} \Big[  2\, L_{\x_1^2}(x)-
L_{\x_1\x_0}(x)\big)\Big]+ {1 \over x+1} \Big[\zeta(2)+
L_{\x_0^2}(x)-2\,L_{\x_0\x_1}(x)\Big]\, dx\ . \end{equation} This is
an integral of functions of weight at most two in a single variable
$x$, in the hyperlogarithm algebra $\MZV\otimes_\Q L(R_4)=\MZV
\otimes_\Q L(\Sigma_x)$. A priori, the regularized limits at
$y=\infty$ which were substituted in at the previous stage could
also have had singularities in $\Sigma_x^+$, i.e., at $x=1$ also,
but the method (the Fubini argument) predicts that the  integrand at the fourth stage
$(\ref{fourthstep})$ cannot, as is indeed the case.

 To complete the calculation, we work in $L(R_4)$. A
further integration  gives:
$$I_G=\Reg_{x=\infty}  \Big(  2\, L_{\x_0\x_1^2}(x)- L_{\x_0\x_1\x_0}(x) + \zeta(2)L_{\x_1}(x)+ L_{\x_1\x_0^2}(x)-2\,L_{\x_1\x_0\x_1}(x)\Big)\ ,$$
Therefore at the fifth and final step, we obtain
$$I_G= 2\,\zeta_{\infty}(\x_0\x_1^2) -\zeta_{\infty}(\x_0\x_1\x_0) + \zeta_{\infty}(\x_1\x_0^2) -2\,\zeta_{\infty}(\x_1\x_0\x_1)\ . $$
Substituting the values given in  $(\ref{Exregzetavalues})$, we
conclude that $$I_G= 6\, \zeta(3)\ .$$
\\

\begin{rem}
One can rewrite $(\ref{fourthstep})$  using only dilogarithms and
logarithms:
\begin{equation}\label{app4step}
I_G=\int_0^\infty {1\over \alpha_4}\Big( -\Li_2(-\alpha_4)
+\log^2(\alpha_4+1) -\log(\alpha_4)\log(\alpha_4+1)  \Big)
\end{equation} $$+{1\over \alpha_4+1}\Big(2\,\Li_2(-\alpha_4)
+{1\over 2} \log^2 (\alpha_4) + \zeta(2) \Big) d\alpha_4\ ,$$ where
each term in the integrand  has singularities contained in $
\{0,-1,\infty\}$. One can obtain $(\ref{app4step})$ directly from
$(\ref{appendixStep3})$ by integrating using the dilogarithm
function, regularizing at infinity, and applying the inversion
relation for $\Li_2$ to arrive at  $(\ref{app4step})$. Although this
gives a substantial shortcut, such a method is ad hoc, and will not
work in a more general setting.


\end{rem}

\end{document}